\documentclass[notitlepage,11pt]{article}
\usepackage{amssymb,amsmath,comment}
\catcode`\@=11
\@addtoreset{equation}{section}

\catcode`\@=12

\usepackage{graphicx}

\def\R{\mathbb{R}}
\def\Rbuco{\R^n\!\setminus\!\{0\}}

\def\S{\mathbb{S}}
\def\f{\varphi}
\def\eps{\varepsilon}

\def\div{{\rm div}}
\def\weak{{\,\rightharpoonup\,}}

\def\p{p}
\def\irn{\int_{\R^n}}

\def\proof{\noindent{\textbf{Proof. }}}
\def\QED{\hfill {$\square$}\goodbreak \medskip}

\newtheorem{Theorem}{Theorem}[section]
\newtheorem{Lemma}[Theorem]{Lemma}
\newtheorem{Proposition}[Theorem]{Proposition}
\newtheorem{Corollary}[Theorem]{Corollary}
\newtheorem{Remark}[Theorem]{Remark}

\begin{document}

\title 
{Weighted Sobolev spaces of radially symmetric functions}

\author{Roberta Musina\footnote{Dipartimento di Matematica ed Informatica, Universit\`a di Udine, via delle Scienze, 206 -- 33100 Udine, Italy. Email: {roberta.musina@uniud.it}.
Partially supported by Miur-PRIN 2009WRJ3W7-001 ``Fenomeni di concentrazione e {pro\-ble\-mi} di analisi geometrica''}}

\date{}

\maketitle

\begin{abstract}
\footnotesize We prove dilation invariant inequalities involving radial functions,
 poliharmonic operators 
and weights that are powers of the distance from the origin.
Then we discuss the existence of extremals and in some cases
we compute the best constants.

\medskip

\medskip

\end{abstract}

\bigskip\bigskip\bigskip

\section{Introduction}
\label{S:Introduction}

The starting point of the present paper is the inequality
\begin{equation}
\label{eq:dis_gen_lin}
\irn|x|^{\alpha}|\nabla^ku|^p~\!dx\ge
c_{\alpha,k,p}~\!
\irn|x|^{\alpha-kp}|u|^p~\!dx
\quad{\forall~u\in C^k_{c,r}(\R^n\!\setminus\!\{0\}).}
\end{equation}
Here $n\ge 2$ and $k\ge 1$ are  integers,  $\alpha\in\R$, $p>1$, 
$C^k_{c,r}(\R^n\!\setminus\!\{0\})$ is the space of radially symmetric functions in $C^k_{c}(\R^n\!\setminus\!\{0\})$,
and
$$
\nabla^k =\begin{cases}
~\Delta^{\!m}&\textrm{if $k=2m$ is even,}\\
\nabla\!\Delta^{\!m}&\textrm{if $k=2m+1$ is odd.}
\end{cases}
$$
Let us briefly describe our main results and our motivations.

First we find out the class of parameters $\alpha, k$ and $p$ such that
 (\ref{eq:dis_gen_lin}) holds with a positive 
best constant $c_{\alpha,k,p}$. If $c_{\alpha,k,p}>0$, then the space
$\mathcal D^{k,\p}_{\rm r}(\R^n;|x|^\alpha dx)$, defined as the
completion of $C^k_{c,r}(\R^n\!\setminus\!\{0\})$ 
with respect to the norm
\begin{equation}
\label{eq:norm_general}
\|u\|_{k,\alpha}=\left(\irn|x|^\alpha|\nabla^k u|^\p~dx\right)^{1/p},
\end{equation}
is continuously embedded in $L^p(\R^n;|x|^{\alpha-kp}dx)$.

Notice that
$\mathcal D^{k,\p}_{\rm r}(\R^n;|x|^\alpha dx)$ is the natural ambient 
space in dealing with radially symmetric solutions to  poliharmonic problems
with weights. Having this application in mind, in the second part of the paper we study
dilation invariant inequalities of the type
\begin{equation}
\label{eq:dis_general}
\irn|x|^{\alpha}|\nabla^ku|^p~\!dx\ge
S_{k,q,j}(\alpha)~\!
\left(\irn|x|^{-\beta_{k-j,q}}|\nabla^ju|^q~\!dx\right)^{p/q}
\end{equation}
for $u\in \mathcal D^{k,\p}_{\rm r}(\R^n;|x|^\alpha dx)$, and we show that
the best constant $S_{k,q,j}(\alpha)$ is positive and achieved. Here 
$j=0,...,k-1$ is an integer,
$\nabla^0 u=u$, $q> p$ is given, and 
\begin{equation}
\label{eq:beta_general}
\beta_{k-j,q}=n-q\frac{n-(k-j)p+\alpha}{p}~\!.
\end{equation}

\bigskip

We remark that (\ref{eq:dis_gen_lin}) is closely related to the double 
weighted Hardy--Littlewood--Sobolev inequality 
(see e.g. Stein and Weiss \cite{SW} and Lieb \cite{Lie}). 
The standard Hardy inequality in \cite{H}, \cite{HLP}
is recovered by choosing $k=1$: it is well known that 
\begin{equation}
\label{eq:H_intro}
\irn|x|^\alpha|\nabla u|^\p~dx\ge\left|\frac{n+\alpha}{p}- 1\right|^\p
\irn|x|^{\alpha-\p}| u|^\p~dx\quad \forall u\in C^1_{c}(\R^n\!\setminus\!\{0\}),
\end{equation}
and that the constant in the right hand side can not be improved.

Second order dilation invariant inequalities have
been largely studied since 1952, when Rellich showed
in \cite{Rel52} (see also \cite{Rel69}) 
that
$$
\int_{\R^n}|\Delta u|^2 dx\ge \left(\frac{n(n-4)}{4}\right)^2
\int_{\R^n} |x|^{-4}|u|^{2 }dx
\quad\forall u\in C^{2}_{c,r}(\R^n\!\setminus\!\{0\})~\!.
$$

In the Hilbertian case $p=2$, Ghoussoub and Moradifam
proved in \cite{GM} that $c_{\alpha,2,2}=
\left|{(n-4+\alpha)(n-\alpha)}/{4}\right|^2$ for any $\alpha\in\R$. We quote also \cite{CM1}, where 
a different approach is used.

Only partial results are available if $p\neq 2$ or $k\ge 3$, see for instance 
Mitidieri \cite{Mit2}, Gazzola-Grunau-Mitidieri
\cite{GGM} (in a non-radial setting) and Adimurthi-Santra \cite{AS}.
Related inequalities
can be found in the above cited papers,  in \cite{AGS}, \cite{CM2}, \cite{CM3},
\cite{dFdSM}, \cite{Mor} and in the references therein.

In the present paper we compute $c_{\alpha,k,p}$ for any
$\alpha, k$ and $p$. We put
\begin{equation}
\label{eq:BGamma}
H_{\!\alpha}=\frac{n+\alpha}{p}-1
~,\quad 
\gamma_{\alpha,h}=\left(\frac{n+\alpha}{p}-h\right)\left(n-2+h-\frac{n+\alpha}{p}\right),
\end{equation}
where $h=2,...,k$ is an integer. Notice that $|H_{\!\alpha}|^p$ is the Hardy constant.
\begin{Theorem}
\label{T:theorem_gen}
Let $\alpha\in\R$, $p>1$ and $k\ge 2$.
The best constant in (\ref{eq:dis_gen_lin}) is given by
\begin{equation}
\label{eq:c}
c_{\alpha,k,p}=\begin{cases}
\quad~~\displaystyle{\prod_{h=1}^m|\gamma_{\alpha,2h} |^{\p}}&\textit{if $k=2m$ is even,}\\
\displaystyle{|H_{\!\alpha}|^p\prod_{h=1}^m|\gamma_{\alpha,2h+1}|^{\p}}&\textit{if $k=2m+1$ is odd.}
\end{cases}
\end{equation}
\end{Theorem}
In case $k=2$
Theorem \ref{T:theorem_gen} implies that
\begin{equation}
\label{eq:R}
\int_{\R^n} |x|^{\alpha}|\Delta u|^\p dx\ge
\left|\left(\frac{n+\alpha}{p}-2\right)\left(n-\frac{n+\alpha}{p}\right)\right|^\p
\int_{\R^n} |x|^{\alpha-2\p}|u|^{\p }dx
\end{equation}
for any $u\in C^2_{c,r}(\R^n\!\setminus\!\{0\})$. 
In fact  (\ref{eq:R}) is a corollary of the next result (with $m=1$)
and of the Hardy inequality. Here we agree that $\Delta^0 u=u$.

\begin{Theorem}
\label{T:radial}
Let $\alpha\in\R$, $\p>1$ and let $m\ge 1$ be a given integer. Then
\begin{equation}
\label{eq:AS}
\int_{\R^n} |x|^{\alpha}|\Delta^m u|^\p dx\ge
\left|~\!n -\frac{n+\alpha}{p}\right|^\p
\int_{\R^n} |x|^{\alpha-\p}|\nabla(\Delta^{m-1}u)|^{\p }dx
\end{equation}
for any $u\in C^{2m}_{c,r}(\R^n\!\setminus\!\{0\})$. The constant in the right hand side is sharp.
\end{Theorem}
Notice that in the singular
case $\alpha=2\p-n$ the best constant in (\ref{eq:R}) vanishes, while
$$
\int_{\R^n} |x|^{2\p-n}|\Delta u|^\p dx\ge
\left|n-2\right|^\p
\int_{\R^n} |x|^{\p-n}|\nabla u|^{\p }dx
$$
for any $u\in C^{2}_{c,r}(\R^n\!\setminus\!\{0\})$.

Theorem \ref{T:radial} and (\ref{eq:H_intro}) provide
 explicit best constants in inequalities of the type
\begin{equation}
\label{eq:tedious}
\irn|x|^{\alpha}|\nabla^ku|^p~\!dx\ge
c_{\alpha,k,j,p}~\!
\irn|x|^{\alpha-(k-j)p}|\nabla^ju|^p~\!dx
\quad{\forall~u\in C^k_{c,r}(\R^n\!\setminus\!\{0\}),}
\end{equation}
for any intermediate case
$j=1,...,k-1$, see Remark \ref{R:tedious}.

Weighted inequalities of order $k\ge 2$ for non radial functions
are more involved. In \cite{GM} and  \cite{CM1}, where $p=2$ and $k=2$
are assumed, it is proved that
$$
\int_{\R^n}|x|^\alpha|\Delta u|^2 dx\ge 
\min_{i\in\mathbb{N}\cup\{0\}}\left|\gamma_{\alpha,2}+i(n-2+i)\right|^{2}
\int_{\R^n} |x|^{\alpha-4}|u|^{2 }dx
$$
for any $u\in C^{2}_{c}(\R^n\setminus\{0\})$. Here 
${\gamma_{\alpha,2}=\left(\frac{n-4+\alpha}{2}\right)\left(\frac{n-\alpha}{2}\right)}$,
accordingly with (\ref{eq:BGamma}). In particular, the best constant vanishes if
$-\gamma_{\alpha,2}$ is an eigenvalue of the Laplace-Beltrami operator
on the sphere. The problem of finding the best constant for weighted Rellich type inequalities
in a non radial setting and for general parameters $\alpha, k, p$ is still open.

Next we direct our attention to semilinear inequalities. Assume $c_{\alpha,k,p}>0$
and take an exponent $q>p$. For any integer $j\in\{0,...,k-1\}$ let 
$S_{k,q,j}(\alpha)$  be the best constant
in (\ref{eq:dis_general}), that is,
\begin{equation}
\label{eq:theorem_semilinear} 
S_{k,q,j}(\alpha)=
\inf_{\scriptstyle u\in\mathcal D^{k,\p}_{\rm r}(\R^n;|x|^\alpha dx)\atop\scriptstyle u\ne 0}
\frac{\displaystyle
\int_{\R^n}|x|^{\alpha}|\nabla^k u|^\p dx}
{\left(\displaystyle\int_{\R^n}|x|^{-\beta_{k-j,q}}| \nabla^j u|^q dx\right)^{\p/q}}~\!.
\end{equation}
In Section
\ref{A:poli} we prove the following existence result (see also Section  \ref{SS:Dk=2} for
shorter proofs in case $k=2$, $j\in\{0,1\}$).

\begin{Theorem}
\label{T:extremalsk}
If $q>p>1$ and $j\in\{0,1,...,k-1\}$, then 
$S_{k,q,j}(\alpha)$
is positive and achieved in $\mathcal D^{k,\p}_{\rm r}(\R^n;|x|^\alpha dx)$.
\end{Theorem}

If $n>kp$, $\alpha=0$, $j=0$ and $q=p^{k*}:=\frac{np}{n-kp}$, then $\beta_{k,q}=0$ and
$S_{k,q,0}(0)$
coincides with the radial Sobolev constant $S_{k,p}^{{k\!*}}$.
Actually,  (\ref{eq:theorem_semilinear}) includes the 
$(k-1)$ best constants $S_{k,p}^*, S_{k,p}^{*\!*},...,S_{k,p}^{(k-1)*}$, that are relative to the embeddings 
$$
\mathcal D^{k,\p}_{\rm r}(\R^n)~\hookrightarrow ~\mathcal D^{j,\frac{np}{n-(k-j)p}}_{\rm r}(\R^n)~,
\quad j=1,..., k-1.
$$
We refer to Remark  \ref{R:newSobolev} for details on this subject.

If $k=1$, $\alpha>p-n$  and
$q>p$, then the infimum $S_{1,q,0}(\alpha)$ is closely related to the
celebrated Caffarelli-Kohn-Nirenbeg inequalities in \cite{CKN}. 
If in addition $p=2$, then the minimizers of
$S_{1,q,0}(\alpha)$ are explicitly known since
the paper \cite{CatWan} by Catrina and Wang
(see also \cite{Lie}). In the next theorem, that will be proved in Section \ref{SS:Dk=1},
we extend the Catrina-Wang uniqueness result to the non Hilbertian case
$p\neq 2$. 

\begin{Theorem}
\label{T:uniqueCKN}
If $\alpha\neq p-n$ and $q>\p$, then  problem
\begin{equation}
\label{eq:ELD1}
-\div(|x|^\alpha|\nabla u|^{\p-2}\nabla u)=|x|^{-n+q\frac{n-\p+\alpha}{\p}}|u|^{q-2}u
\quad\textrm{on $\R^n$}
\end{equation}
has a  nontrivial radial
 solution $U$ in $\mathcal D^{1,\p}_{\rm r}(\R^n;|x|^\alpha dx)$ which is unique up to a change
 of sign and up to a transform of the type
$u(x)\mapsto \rho^{\frac{n-p+\alpha}{p}}~\!u(\rho x)$,  $\rho>0$.
Moreover,
$U$ achieves the best constant $S_{1,q,0}(\alpha)$
and it is given by
$$
U(|x|)= \left(\frac{q}{p}~\!
\frac{\left|{n-p+\alpha}\right|^p}{(p-1)^{p-1}}\right)^{\frac{1}{q-\p}}
\left(1+|x|^{\frac{(n-\p+\alpha)(q-\p)}{\p(\p-1)}}\right)^{\frac{\p}{\p-q}}.
$$
\end{Theorem}
Notice that no a-priori assumption on the sign of $u$ is needed in Theorem \ref{T:uniqueCKN}.

Now assume  $k=2$ and  
$\alpha\notin\{2\p-n,n\p-n\}$. Then $c_{\alpha,2,p}>0$ by Theorem \ref{T:theorem_gen}
and for any $q>p$ the infima $S_{2,q,0}(\alpha)$, $S_{2,q,1}(\alpha)$
are both positive and achieved by Theorem \ref{T:extremalsk}. In Section \ref{SS:Dk=2} we prove our main result concerning the best constant
\begin{equation}
\label{eq:S2q1}
S_{2,q,1}(\alpha)=
\inf_{\scriptstyle u\in\mathcal D^{2,\p}_{\rm r}(\R^n;|x|^\alpha dx)\atop\scriptstyle u\ne 0}
\frac{\displaystyle
\int_{\R^n}|x|^{\alpha}|\Delta u|^\p dx}
{\left(\displaystyle\int_{\R^n}|x|^{-\beta_{1,q}}| \nabla u|^q dx\right)^{\p/q}}~\!,
\end{equation}
where $\beta_{1,q}=n-q\frac{n-p+\alpha}{p}$, accordingly with (\ref{eq:beta_general}).
We identify functions that coincide up to 
a change of sign and a transform of the type
\begin{equation}
\label{eq:new}    
u(x)\mapsto \rho^{\frac{n-2p+\alpha}{p}}~\!u(\rho x)~,~~\rho>0.
\end{equation}


\begin{Theorem}
\label{T:RSachievedJ}
If $\alpha\notin\{2p-n,np-n\}$ and $q>p$, then  the equation
\begin{equation}
\label{eq:gianni}
\Delta\left(|x|^\alpha|\Delta U|^{p-2}\Delta U\right)+\div\left(|x|^{-\beta_{1,q}}|\nabla U|^{q-2}\nabla U\right)
=0
\end{equation}
has a unique 
nontrivial  solution $U\in \mathcal D^{2,\p}_{\rm r}(\R^n;|x|^{\alpha} dx)$.
More precisely, $U$ achieves the best constant $S_{2,q,1}(\alpha)$ and it is given by
\begin{eqnarray*}
&&U(x)=k~\displaystyle\int_{|x|}^\infty t^{1-\frac{\alpha}{\p-1}}
\left(1+t^{\frac{(n\p-n-\alpha)(q-\p)}{\p(\p-1)}}\right)^{\!\frac{\p}{\p-q}}~\!dt\quad
\textit{if $\alpha>2\p-n$}\\
&&U(x)=k~\!
\displaystyle\int_{0}^{|x|} t^{1-\frac{\alpha}{\p-1}}
\left(1+t^{\frac{(n\p-n-\alpha)(q-\p)}{\p(\p-1)}}\right)^{\frac{\p}{\p-q}}~\!dt\quad
\textit{if $\alpha<2\p-n$,}
\end{eqnarray*}
where
$$
k=\left(\frac{q}{p}~\!\frac{\left|{np-n-\alpha}\right|^p}{(\p-1)^{\p-1}}\right)^{\frac{1}{q-\p}}
~\!.
$$
\end{Theorem}
Thanks to Theorem \ref{T:RSachievedJ}, the sharp value of $S_{2,q,1}(\alpha)$
is  computable in terms of Gamma functions, as in
\cite{Au}, \cite{CatWan} and \cite{Ta}.  
We refer to Remarks \ref{R:degenerate1} and \ref{R:degenerate2} for few
comments about the "critical cases" $\alpha=2p-n$, $\alpha=np-n$.

Let us point out a corollary of Theorem \ref{T:RSachievedJ}
that deals with the limiting embedding $\mathcal D^{2,p}(\R^n)\hookrightarrow
\mathcal D^{1,p^*}(\R^n)$, where $p>2n$ and $p^*=\frac{np}{n-p}$ is the (first order)
critical exponent (see also Remark \ref{R:critico} in Section \ref{SS:Dk=2}).

As usual, we denote by 
$\Delta_{p^*}=\div(|\nabla\cdot|^{p^*-2}\nabla\cdot)$ the $p^*$-Laplace operator. 
\begin{Corollary}
\label{C:Malchiod}
Assume that $n> 2p$. Then the problem
\begin{equation}
\label{eq:ELSobolev}
\Delta\left(|\Delta U|^{p-2}\Delta U\right)+\Delta_{p^*} U=0
\end{equation}
has a unique nontrivial 
solution $U\in \mathcal D^{2,p}_{\rm r}(\R^n)$. More precisely, $U$ achieves
\begin{equation}
\label{eq:DN_alcritico}
{S_{2,p}^*}:=\inf_{\scriptstyle 
u\in \mathcal D^{2,p}_{\rm r}(\R^n)\atop\scriptstyle u\ne 0}
\frac{\displaystyle
\int_{\R^n}\left|\Delta u\right|^{p} dx}
{\left(\displaystyle\int_{\R^n}|\nabla u|^{p^*}~\!dx\right)^{p/p^*}}
\end{equation}
and it is given by
$$
U(x)= n^{\frac{n-p}{p}}~\!\left(\frac{n(p-1)}{n-p}\right)^{\frac{n-p}{p^2}}
\int_{|x|}^{\infty}s\left(1+s^{p^*}\right)^{\!\!\frac{p-n}{p}}~\!ds.
$$
\end{Corollary}

\bigskip

There are a few ways to prove inequalities like (\ref{eq:dis_gen_lin}) or (\ref{eq:dis_general}).
In \cite {Mit2} Mitidieri applied his powerful {\em "simple approach"} to compute $c_{\alpha,2,p}$, among other best constants,
when $\gamma_{\alpha,2}\ge 0$. Pointwise estimates 
 are frequently used to obtain integral inequalities, see for instance the papers \cite{CKN} by
Caffarelli, Kohn and Nirenberg and the more recent \cite{AGS}, \cite{AS}. In presence
of symmetries, Calanchi-Ruf  \cite{CR} and de Figueiredo-dos Santos-Miyagaki  \cite{dFdSM} obtained embedding
results as corollaries of a  {\em radial lemma} 
 (in the spirit of \cite{PLL0} and \cite{Ni}).

 Here we use a different approach. We start in Section \ref{S:simple} by proving
 Theorems  \ref{T:radial} and \ref{T:theorem_gen} via the Hardy inequality
for functions of one real variable. No pointwise estimates are needed. 
To study  (\ref{eq:dis_general}) we argue in the opposite direction
with respect to the above mentioned papers:
first we prove certain embedding theorems, and then we  infer
the desired inequalities. We focus our attention on (\ref{eq:dis_general}) even if we can obtain also 
pointwise estimates, radial lemmas and informations on the regularity of
radial functions, see Remark \ref{R:regularity}. Roughly speaking, to get more inequalities in case
$c_{\alpha,k,p}> 0$ we define in Section \ref{S:Rellich}
a $k$-th order Emden-Fowler transform
$\mathcal D^{k,\p}_{\rm r}(\R^n;|x|^\alpha dx)\to W^{k,p}(\R)$ and we show that
the induced norm on $W^{k,p}(\R)$
is equivalent to the standard one. Then,  classical results about the space $W^{k,p}(\R)$
provide embedding theorems for $\mathcal D^{k,\p}_{\rm r}(\R^n;|x|^\alpha dx)$, and the
inequalities we are interested in readily follow.

The first step in this program consists in highlighting a suitable class of equivalent
norms on the Sobolev spaces $W^{k,\p}(\R)$. We start with the
lowest indexes
$k=1$ and $k=2$ in Sections \ref{S:Preliminaries} and \ref{S:W2p}, respectively.
The higher order case $k\ge 3$ will be briefly discussed in Section \ref{S:higher}.

\bigskip

\small
\noindent
{\bf Notation.}~\! 
We denote by $c$ any nonnegative universal constant. 

We set $\R_+=(0,\infty)$. For any integer $n\ge 2$ we denote by 
$\omega_n$ the $n-1$ dimensional measure of the unit sphere
$\S^{n-1}$ in $\R^n$.

If $\Omega\subseteq\R^n$ is a rotationally invariant domain and $k\ge 0$ is an integer, we denote by
$C^k_{c,r}(\Omega)$ the space of radially symmetric functions $u\in C^k_c(\Omega)$. 
For $u\in C^1_{c,r}(\Omega)$ we let $u'$ be
the radial derivative of $u$. Thus $|x|u'(x)=\nabla u(x)\cdot x$.

The exponent $q'=\frac{q}{q-1}$ is conjugate exponent to $q\in(1,\infty)$.

Let  
$\omega$ be a non-negative measurable function on 
a domain $\Omega\in\R^n$, $n\ge 1$. The weighted
Lebesgue space $L^q(\Omega;\omega(x)~\! dx)$ is the space of measurable maps
$u$ in $\Omega$ with finite norm $\left(\int_{\Omega}|u|^q
\omega(x)~\!dx\right)^{1/q}$.  For $\omega\equiv 1$ we denote by 
$\|u\|_q$
the standard norm in
$L^q(\Omega)=L^q(\Omega;~\! dx)$. 

The
norm in the Sobolev space $W^{k,q}(\R)$ is given by
$$
\|g\|_{W^{k,q}}=\left(\int_\R|g^{{\rm (k)}}|^q~\!ds+\int_\R|g|^q~\!ds\right)^{1/q}.
$$
If $n>kp$ then the space 
 $\mathcal D^{k,\p}(\R^n)$ is the closure of $C^\infty_c(\R^n)$
 with respect to the norm
 $$
 \|u\|=\left(\irn|\nabla^k u|^\p~dx\right)^{1/p}.
 $$
 We put 
 $\mathcal D^{k,\p}_r(\R^n)=
\left\{u\in \mathcal D^{k,\p}(\R^n)~|~u=u(|x|)\right\}$.

 Let $\p^{{k\!*}}=\frac{np}{n-k\p}$ be the $k$-th order critical exponent.
 The radial Sobolev constant
 $$
 S_{k,p}^{{k\!*}}:=\inf_{\scriptstyle 
u\in \mathcal D^{k,p}_r(\R^n)\atop\scriptstyle u\ne 0}
\frac{\displaystyle
\int_{\R^n}\left|\nabla^k u\right|^{p} dx}
{\left(\displaystyle\int_{\R^n}| u|^{p^{{k\!*}}}~\!dx\right)^{p/p^{{k\!*}}}}
$$
is positive and achieved (see also Theorem \ref{T:extremalsk} and Remark
\ref{R:newSobolev} in Section \ref{SS:Dk=2} below).

Assume that $k=1$ or $p=2$. Then it is well known that 
$S_{k,p}^{{k\!*}}$ is the best constant in the embedding
$\mathcal D^{k,\p}(\R^n)\hookrightarrow L^{\p^{{k\!*}}}(\R^n)$, that is,
$$
S_{k,p}^{{k\!*}}=\inf_
{\scriptstyle u\in \mathcal D^{k,p}(\R^n)\atop\scriptstyle u\ne 0}
\frac{\displaystyle
\int_{\R^n}\left|\nabla^k u\right|^{p} dx}
{\left(\displaystyle\int_{\R^n}| u|^{p^{{k\!*}}}~\!dx\right)^{p/p^{{k\!*}}}}~\!,
$$
see \cite{Au}, \cite{Ta} and for instance \cite{GGS} for the poliharmonic case.
The notation $k*$ means $\underbrace{**...*}_{k~\!\textrm{times}}$. 
If $k\in\{1,2\}$ we write $*, {*\!*}$ instead
 of ${1*}, {2*}$, respectively.

\normalsize

\section{Higher order Hardy-Rellich inequalities}
\label{S:simple}

Throughout this paper we will use the Hardy inequality for  functions in $C^1_c(\R_+)$ several times.
We recall that for any $a\in\R$, $\p>1$ and $\omega\in C^1_c(\R_+)$ the inequality
\begin{equation}
\label{eq:H}
\int_0^\infty r^a|\omega'|^\p dr\ge \left|\frac{a+1-\p}{\p}\right|^\p
\int_0^\infty r^{a-\p}|\omega|^\p dr
\end{equation}
holds with a sharp and non achieved constant in the right hand side, see \cite{H}, \cite{HLP}. 
We notice that (\ref{eq:H}) holds as well in the "singular case"
$a=p-1$, with null best constant. For, use
$\omega_\eps({r})=\omega(r^\eps)$  as test function, where
$\omega\in C^1_c(\R_+)\setminus\{0\}$ is fixed, to get
$$
\frac{\displaystyle\int_0^\infty r^{p-1}|\omega_\eps'|^\p dr}{\displaystyle\int_0^\infty r^{-1}|\omega_\eps|^\p dr}=
\eps^p~ \frac{\displaystyle\int_0^\infty r^{p-1}|\omega'|^\p dr}{\displaystyle\int_0^\infty r^{-1}|\omega|^\p dr}
\longrightarrow 0\quad\textrm{as $\eps\to 0$.}
$$

We point out a simple but very useful  corollary to the
Hardy inequality.

\begin{Lemma}
\label{L:1dim}
Let $\tau, {\lambda}\in\R$, $\p>1$ and   $v\in C^2_c(\R_+)$. Then the inequalities
\begin{gather}
\label{eq:i)}
\int_0^\infty r^\tau\left|v''+({\lambda}-1)r^{-1}v'\right|^\p~\!dr
\ge
\left|\frac{\tau+1-{\lambda}\p}{\p}\right|^\p
\int_0^\infty r^{\tau-\p}\left|v'\right|^\p~\!dr\\
\label{eq:ii)}
\int_0^\infty\!\!\! r^\tau\left|v''+({\lambda}-1)r^{-1}v'\right|^\p dr
\ge
\left|\frac{(\tau+1-{\lambda}\p)(\tau+1-2\p)}{\p^2}\right|^\p
\int_0^\infty \!\!\! r^{\tau-2\p}|v|^\p dr
\end{gather}
hold with sharp constants.
\end{Lemma}

\proof
Inequalities (\ref{eq:i)}), (\ref{eq:ii)}) are immediate consequences of (\ref{eq:H}), since for any $v\in C^2_c(\R_+)$ we have that
\begin{eqnarray*}
\int_0^\infty r^\tau\left|v''+({\lambda}-1)r^{-1}v'\right|^\p dr&=&
\int_0^\infty r^{\tau-(\lambda-1)p}\left|\left(r^{\lambda-1}v'\right)'\right|^\p dr
\\
&\ge&
\left|\frac{\tau+1-\lambda p}{p}\right|^p
\int_0^\infty r^{\tau-\p}\left|v'\right|^\p~\!dr \\
&\ge&
\left|\frac{\tau+1-\lambda p}{p}\right|^p\left|\frac{\tau+1-2 p}{p}\right|^p
\!\!
\int_0^\infty r^{\tau-2\p}\left|v\right|^\p~\!dr~\!.
\end{eqnarray*}
To prove the sharpness of the constants in (\ref{eq:i)}), (\ref{eq:ii)}) 
we fix a nontrivial function
$v\in C^2_c(\R_+)$ and we put 
$$v_\eps({r})=\eps^{\frac{1}{p}}r^{2-\frac{1+\tau}{p}}v(r^\eps),$$ where
$\eps\to 0^+$. To conclude,
it suffices to compute
\begin{eqnarray*}
\int_0^\infty\!\!\! r^\tau\left|v_\eps''+({\lambda}-1)r^{-1}v_\eps'\right|^\p~\!dr&=&
\left|\frac{(\tau+1-{\lambda}\p)(\tau+1-2\p)}{\p^2}\right|^\p
\int_0^\infty\!\!\! r^{-1}|v|^p~\!dr+o(1)\\
\int_0^\infty r^{\tau-p}|v'_\eps|^p~\!dr&=&
\left|\frac{\tau+1-2\p}{\p}\right|^\p \int_0^\infty r^{-1}|v|^p~\!dr+o(1)\\
\int_0^\infty \!\!\! r^{\tau-2\p}|v_\eps|^\p dr&=&
\int_0^\infty r^{-1}|v|^p~\!dr.
\end{eqnarray*}
\QED

\subsection{Proof of Theorem \ref{T:radial}}
Fix any
$u\in C^{2m}_{c,r}(\R^n\setminus\{0\})$. Since
\begin{gather*}
\int_{\R^n} |x|^{\alpha}|\Delta u|^\p dx=
\omega_n
\int_0^\infty r^{n-1+\alpha}|u''+(n-1)r^{-1}u'|^\p dr\\
\int_{\R^n} |x|^{\alpha-\p}|\nabla u|^{\p }dx=\omega_n
\int_0^\infty r^{n-1+\alpha-\p}|u'|^\p dr,
\end{gather*}
then from Lemma \ref{L:1dim} it follows that
\begin{equation}
\label{eq:AS1}
\int_{\R^n} |x|^{\alpha}|\Delta u|^\p dx\ge
\left|~\!n -\frac{n+\alpha}{\p}\right|^\p
\int_{\R^n} |x|^{\alpha-\p}|\nabla u|^{\p }dx
\end{equation}
and that the constant in the right hand side is the best possible.
Thus (\ref{eq:AS}) is proved when $m=1$. 

To conclude the proof in case
$m>1$ it suffices 
to write down (\ref{eq:AS1})
with $\Delta^{m-1}u\in C^{2}_{c,r}(\R^n\!\setminus\!\{0\})$ instead of $u$.
\QED

\subsection{Proof of Theorem \ref{T:theorem_gen}}
Let $\tilde c_{\alpha,k,p}$ be the best constant in (\ref{eq:dis_gen_lin}), and
let $c_{\alpha,k,p}$ be the constant defined in  (\ref{eq:c}). 
We start by proving that $\tilde c_{\alpha,k,p}\ge c_{\alpha,k,p}$ in case
$k=2m$ is an even integer. We have to show that
\begin{equation}
\label{eq:indu_pari}
\irn|x|^{\alpha}|\Delta^m u|^p~\!dx\ge
\left(\prod_{h=1}^m|\gamma_{\alpha,2h} |^{\p}\right)
\irn|x|^{\alpha-2mp}|u|^p~\!dx
\quad{\forall~u\in C^{2m}_{c,r}(\R^n\!\setminus\!\{0\}).}
\end{equation}
If $m=1$ then (\ref{eq:indu_pari}) reduces to (\ref{eq:R}), that is an
immediate consequence of (\ref{eq:AS1}) and of the Hardy inequality
(\ref{eq:H_intro}). 
Assume that (\ref{eq:indu_pari}) holds for some $m\ge1$
and for any $\alpha\in\R$. Fix
$u\in C^{2m+2}_{c,r}(\R^n\!\setminus\!\{0\})$ and use (\ref{eq:R}) to infer
$$
\irn|x|^{\alpha-2mp}|\Delta u|^p~\!dx\ge
|\gamma_{\alpha,2m+2} |^{\p}~\!
\irn|x|^{\alpha-2(m+1)p}|u|^p~\!dx,
$$
since $\gamma_{\alpha-2mp,2}=\gamma_{\alpha,2m+2}$.
Thus we have that
\begin{eqnarray*}
\irn|x|^{\alpha}|\Delta^{m+1} u|^p~\!dx&=&
\irn|x|^{\alpha}|\Delta^m(\Delta u)|^p~\!dx\\
&\ge& \left(\prod_{h=1}^m|\gamma_{\alpha,2h} |^{\p}\right)
\irn|x|^{\alpha-2mp}|\Delta u|^p~\!dx\\
&\ge& \left(\prod_{h=1}^m|\gamma_{\alpha,2h}|^{\p}\right) |\gamma_{\alpha,2m+2} |^{\p}~\!
\irn|x|^{\alpha-2(m+1)p}|u|^p~\!dx,
\end{eqnarray*}
as desired. If $k=2m+1$ is odd we use the Hardy inequality
and the first part of the proof to get
\begin{eqnarray*}
\irn|x|^{\alpha}|\nabla(\Delta^m u)|^p~\!dx&\ge&
|H_{\!\alpha}|^p
\irn|x|^{\alpha-\p}|\Delta^m u|^\p~dx\\
&\ge&|H_{\!\alpha}|^p \left(\prod_{h=1}^m|\gamma_{\alpha-p,2h}|^{\p}\right) 
\irn|x|^{\alpha-\p}|\Delta^m u|^\p~dx
\end{eqnarray*}
for any $u\in C^{2m+1}_{c,r}(\R^n\!\setminus\!\{0\})$.
Thus $\tilde c_{\alpha,k,p}\ge c_{\alpha,k,p}$ also in this case, since
 $\gamma_{\alpha-p,2h}=
\gamma_{\alpha,2h+1}$ for any integer $h$. 

To prove that $\tilde c_{\alpha,k,p}\le c_{\alpha,k,p}$,  fix a nontrivial function
$u\in C^k_{c,r}(\R^n\!\setminus\!\{0\})$. For any $\eps>0$  define the radial function
$u_\eps(|x|)=|x|^{k-\frac{n+\alpha}{p}}u(|x|^\eps)$. Direct computations and induction can be used to check that
\begin{eqnarray*}
\irn|x|^{\alpha}|\nabla^ku_\eps|^p~\!dx&=&{c_{\alpha,k,p}}~\!{\eps}^{-1}
\irn |x|^{-n}\left|u+\eps\psi_\eps\right|^p~\!dx\\
\irn|x|^{\alpha-kp}|u_\eps|^p~\!dx&=&{\eps}^{-1}
\irn|x|^{-n}|u|^p~\!dx,
\end{eqnarray*}
where the radial function $\psi_\eps\in C^0_c(\R^n)$ is a linear combination of
derivatives of $u$, such that $\sup_\eps\|\psi_\eps\|_\infty<\infty$.
The conclusion readily follows.
\QED

\begin{Remark}
\label{R:alternative0}
An alternative proof of Theorem \ref{T:theorem_gen} is suggested in 
Remark \ref{R:alternative}. 
\end{Remark}

\begin{Remark}
\label{R:tedious}
Similar computations allow us to find the best constant $c_{\alpha,k,j,p}$ in (\ref{eq:tedious}).
Theorem \ref{T:radial} and the Hardy inequality provide the values of
$c_{\alpha,k,j,p}$ for $j=k-1$ and $j=k-2$. For smaller indexes we put
$\delta_{j}=n-1-\frac{n+\alpha}{p}+k-j$.
If $k=2m$ is even we have 
$$
c_{\alpha,2m,j,p}=\left\{
\begin{array}{cl}
\displaystyle{\prod_{h=1}^{m-i}|\gamma_{\alpha,2h} |^{\p}}&\textit{if $j=2i$}\\
\displaystyle{|\delta_{j}|^p}
\displaystyle{\prod_{h=1}^{m-i-1}|\gamma_{\alpha,2h}|^{\p}}&\textit{if $j=2i+1$,}
\end{array}
\right. 
$$
while if $k=2m+1$ it results that
$$
c_{\alpha,(2m+1),j,p}=\left\{
\begin{array}{cl}
\displaystyle{|H_\alpha|^p}
\displaystyle{\prod_{h=1}^{m-i}|\gamma_{\alpha,2h+1}|^{\p}}&\textit{if $j=2i$,}\\
\displaystyle{|H_\alpha\delta_{j}|^p}
\displaystyle{\prod_{h=1}^{m-i-1}|\gamma_{2h+1}|^{\p}}&\textit{if $j=2i+1$.}
\end{array}
\right. 
$$
\end{Remark}

\section{The Emden-Fowler transform\\ and first order inequalities}
\label{S:Rellich}

From now on we will assume that $c_{\alpha,k,p}$ is positive, that is,
\begin{equation}
\label{eq:nondeg_gen}
\begin{cases}
\gamma_{\alpha,2h}\neq 0\quad\quad\forall h=1,...,m&\textrm{if $k=2m$ is even,}\\
\gamma_{\alpha,2h+1}\neq 0\quad\forall h=1,...,m~~
\textrm{and}~~H_{\!\alpha}\neq 0&\textrm{if $k=2m+1$ is odd.}
\end{cases}
\end{equation}
We study the properties of the Banach space
$\mathcal D^{k,\p}_{\rm r}(\R^n;|x|^\alpha dx)$ by using a suitable transform
$W^{k,p}(\R)\to \mathcal D^{k,\p}_{\rm r}(\R^n;|x|^\alpha dx)$ and by taking
advantage of the results in the previous section. More precisely, for any
$k\ge 0$
we define the (inverse)
$k$-th order Emden-Fowler transform
\begin{equation}
\label{eq:EF}
\mathcal T_k:C^k_c(\R)\to C^k_{c,r}(\Rbuco)~\quad
(\mathcal T_k g)(x)=|x|^{-{H}_{\alpha,k}}g(-\log|x|),
\end{equation}
where we have set
\begin{equation}
\label{eq:B}
{H}_{\alpha,k}= \frac{n+\alpha}{\p}-k.
\end{equation}
Notice that ${H}_{\alpha,1}=H_\alpha$, compare with (\ref{eq:BGamma}).
Then we show that $\mathcal T_k$ extends 
to a bicontinuous isomorphism
$W^{k,p}(\R)\to \mathcal D^{k,\p}_{\rm r}(\R^n;|x|^\alpha dx)$ for any
$k\ge 1$.  A crucial step in this program consists in finding out a
large class of equivalent norms in $W^{k,p}(\R)$.
We start with the lower order case $k=1$. The cases $k=2$ and
$k\ge 3$ will be studied in the last two sections.

\subsection{Equivalent norms on $W^{1,\p}(\R)$}
\label{S:Preliminaries}
We point out a simple lemma, based on the Hardy inequality for functions in $C^1_c(\R_+)$.

\begin{Lemma}
\label{L:1dim1st}
Let $p>1$ and ${\lambda}\in\R$. Then
$$
M_p(\lambda):=\inf_{\scriptstyle f\in W^{1,\p}(\R)\atop\scriptstyle f\ne 0}
\frac{\displaystyle
\int_\R\left|f'-{\lambda}~\!f~\!\right|^\p ds}
{\displaystyle\int_\R|f|^\p~\!ds}=|\lambda|^\p.
$$
\end{Lemma}

\proof
To any $f\in C^1_c(\R)$ we
associate the function $v(r):=r^{{\lambda}}f(-\log r)$. Then
$v\in C^1_c(\R_+)$ and a direct computation shows that
$$
\frac{\displaystyle
\int_\R\left|f'-{\lambda}~\!f~\!\right|^\p ds}
{\displaystyle\int_\R|f|^\p~\!ds}=
\frac{\displaystyle
\int_0^\infty r^{(1-{\lambda})\p-1}|v'|^\p dr}
{\displaystyle\int_0^\infty r^{-{\lambda}\p-1}|v|^\p~\!dr}~\!.
$$
The conclusion readily follows by using (\ref{eq:H}) and a density argument.
\QED
Now we take an exponent $q>p$ and we study the infimum
\begin{equation}
\label{eq:M}
M_{\p,q}(\lambda):=\inf_{\scriptstyle f\in W^{1,\p}(\R)\atop\scriptstyle f\ne 0}
\frac{\displaystyle
\int_\R\left|f'-{\lambda}~\!f~\!\right|^\p ds}
{\left(\displaystyle\int_\R|f|^q~\!ds\right)^{\p/q}}~\!.
\end{equation}

\begin{Remark}
\label{R:M=0}
A standard rescaling argument can be used to check that
$M_{\p,q}(0)=0$. 
\end{Remark}
The next proposition gives us the
equivalent norms we need in case $k=1$. Its proof is immediate,
by Lemma \ref{L:1dim1st} and by Sobolev embedding theorem.

\begin{Proposition}
\label{P:equivalent1}
Let $\p>1$ and $\lambda\in\R\setminus\{0\}$. Then 
$$
\|f\|:= 
\left(\int_\R|f'-\lambda f|^\p~\!ds\right)^{1/\p}
$$
is equivalent to the standard norm on $W^{1,\p}(\R)$. Thus, for any
$q> \p$ the infimum $M_{\p,q}(\lambda)$ is positive.
\end{Proposition}
\begin{Remark}
\label{R:standard}
The minimization problem in (\ref{eq:M}) is non compact, due to translations
in $\R$. By
nowadays standard arguments one can prove that for
every bounded minimizing sequence $f_h$, there exists
a sequence $s_h$ in $\R$ such that $f_h(\cdot-s_h)$ is relatively compact
in $W^{1,p}(\R)$. Hence, 
$M_{\p,q}(\lambda)$ is attained by some function $f\neq 0$ which solves
\begin{equation}
\label{eq:EL1}
-\left(|f'-\lambda f|^{\p-2}(f'-\lambda f)\right)'-\lambda|f'-\lambda f|^{\p-2}(f'-\lambda f)=
|f|^{q-2}f
\quad\textrm{on $\R$}
\end{equation}
up to a Lagrange multiplier.
\end{Remark}

If $\p=2<q$, then extremals for
$$
M_{2,q}(\lambda)=\inf_{\scriptstyle f\in H^{1}(\R)\atop\scriptstyle f\ne 0}
\frac{\displaystyle
\int_\R\left|f'-{\lambda}~\!f~\!\right|^2 ds}
{\left(\displaystyle\int_\R|f|^q~\!ds\right)^{2/q}}=
\inf_{\scriptstyle f\in H^1(\R)\atop\scriptstyle f\ne 0}
\frac{\displaystyle
\int_\R\left(|f'|^2+{\lambda}^2~\!|f|^2~\!\right) ds}
{\left(\displaystyle\int_\R|f|^q~\!ds\right)^{2/q}}
$$
give rise to nontrivial solutions of the Emden-Fowler (or Schr\"odinger) equation
\begin{equation}
\label{eq:EF2}
-f''+\lambda^2 f=|f|^{q-2}f\quad\textrm{on $\R$.}
\end{equation}
It has been shown in \cite{CatWan} that, up to translations, equation (\ref{eq:EF2}) has a unique
positive solution $F\in H^1(\R)$, which is explicitly known. The interest of Catrina and Wang in the ODE (\ref{eq:EF2}) was motivated by
its relevance with the Caffarelli-Kohn-Nirenberg inequalities in the Hilbertian case $p=2$.

Now we state a uniqueness result for nontrivial solutions
$f\in  W^{1,p}(\R)$
to (\ref{eq:EL1}). Notice that we
do not require any sign assumption on $f$. Thus
some care is needed, as the
exponents $p$,  $q$ might be smaller than $2$.

We identify functions that coincide
up to a translation and a change of sign.

\begin{Theorem}
\label{T:EF1}
Let $q>\p>1$ and $\lambda\in\R\setminus\{0\}$. Then the ordinary
differential equation
(\ref{eq:EL1}) has a unique nontrivial solution $F\in W^{1,\p}(\R)$.
More precisely, $F$ achieves the best constant
$M_{\p,q}(\lambda)$, and it is given by
\begin{equation}
\label{eq:numerare}
F(s)=\left(q\left(\frac{\p}{\p-1}\right)^{\p-1}\left|\frac{\lambda}{2}\right|^\p\right)^{\frac{1}{q-\p}}
~\!e^{\frac{\lambda(\p-2)}{2(\p-1)}s}~\!
\left(\cosh\left(\frac{\lambda(q-\p)}{2(\p-1)}~\!s\right)~\!\right)^{\frac{\p}{\p-q}}.
\end{equation}
\end{Theorem}

\proof
Let $f\in W^{1,\p}(\R)\setminus\{0\}$ be a solution to (\ref{eq:EL1}),
and put 
$$
\f= |f'-\lambda f|^{\p-2}(f'-\lambda f).
$$
The pair $f,\f$ solves
\begin{equation}
\label{eq:ODE1}
\begin{cases}
f'-\lambda f=|\f|^{\p'-2}\f\\
-\f'-\lambda\f=|f|^{q-2}f
\end{cases}
\end{equation}
in the sense of distributions.
Notice that   $\p',q$ satisfy the standard anticoercivity assumption $(\p'-1)(q-1)>1$. 
Clearly,   $\f\in L^{\p'}(\R)$ and
$\|\f\|_{\p'}=\|f'-\lambda f\|_{\p}$. Since $-\f'=\lambda\f+|f|^{q-2}f\in L^{p'}(\R)$
by Sobolev embeddings, then $\f\in W^{1,p'}(\R)$. Thus
$\f\in C^1(\R)$, as $f$ and $\f$ are continuous function. But then
also $f$ is of class $C^1$, since $f'=\lambda f+|\f|^{\p'-2}\f$.
Thus the pair $f,\f$ is a classical homoclinic solution to (\ref{eq:ODE1}).

The system  (\ref{eq:ODE1}) is conservative, with with Hamiltonian energy
$$H(f,\f)=\lambda f\f+\frac{1}{q}|f|^q+\frac{1}{\p'}|\f|^{\p'}.$$
In particular, (\ref{eq:ODE1}) is equivalent to
\begin{equation}
\label{eq:HamSys}
\begin{cases}
f'=\partial_{2}H(f,\f)\\
\f'=-\partial_{1}H(f,\f).
\end{cases}
\end{equation}
From $f\in  W^{1,\p}(\R), \f\in W^{1,p'}(\R)$ one infers that $f,\f$ vanish at infinity, and therefore
\begin{equation}
\label{eq:Ham1}
\lambda f\f+\frac{1}{q}|f|^q+\frac{1}{\p'}|\f|^{\p'}=0.
\end{equation}
Notice that $\lambda f\f<0$ on the set $\{f\neq 0\}=\{\f\neq 0\}$.
We can assume that $f$ achieves its positive maximum at some point $s_0$. Using $f'(s_0)=0$,
(\ref{eq:ODE1}) and (\ref{eq:Ham1}) one can uniquely compute
the values of $f(s_0)> 0$ and $\lambda\f(s_0)< 0$. 
Since for any initial datum $f_0> 0$, $\f_0\neq 0$
the Cauchy problem for 
(\ref{eq:HamSys}) has a unique local
solution, to conclude the proof
we only have to show that 
the pair
$F,\Phi$ solves (\ref{eq:ODE1}), where 
$$
\Phi=|F'-\lambda F|^{\p-2}(F'-\lambda F).
$$
In order to avoid long computations one can argue as follows. Put
$$
k=\left(q\left(\frac{\p}{\p-1}\right)^{\p-1}\left|\frac{\lambda}{2}\right|^\p\right)^{\frac{1}{q-\p}}~,\quad
c_1=\frac{\lambda(\p-2)}{2(\p-1)}~,\quad c_2=\frac{\lambda(q-\p)}{2(\p-1)},
$$
so that
$
F(s)=k
~\!e^{c_1s}~\!
\left(\cosh c_2~\!s~\!\right)^{\frac{\p}{\p-q}}
$, and compute
$$
\Phi=-\left(k~\!\frac{\p}{2(p-1)}\right)^{\p-1}|\lambda|^{\p-2}\lambda
~\!e^{(\p-1)(c_1+c_2)s}~\!
\left(\cosh c_2~\!s~\!\right)^{\frac{q(\p-1)}{\p-q}}.
$$
Now it is easy to check that the pair $F,\Phi$ satisfies the conservation law (\ref{eq:Ham1}),
that is sufficient to conclude that $F,\Phi$ solves (\ref{eq:ODE1}), as desired. 
\QED

\subsection{The space $\mathcal D^{1,\p}_{\rm r}(\R^n;|x|^\alpha dx)$}
\label{SS:Dk=1}
If $\alpha\neq p-n$, then the Banach space 
$\mathcal D^{1,\p}_{\rm r}(\R^n;|x|^\alpha dx)$, endowed with the norm
$$
\|u\|^p_{1,\alpha}=\irn|x|^\alpha|\nabla u|^\p~dx,
$$
is continuously embedded into $L^p(\R^n;|x|^{\alpha-p}dx)$
by the Hardy inequality.

For $g\in C^1_c(\R)$ we put $(\mathcal T_1 g)(x)=|x|^{-H_{\!\alpha}}g(-\log|x|)$, where
$H_{\!\alpha}= \frac{n+\alpha}{p}-1$.
Then clearly $\mathcal T_1:  C^1_c(\R)\to C^1_{c,r}(\Rbuco)$ is a linear, invertible transform.
\begin{Lemma}
\label{L:space_explicit1}
Assume $\alpha\neq p-n$.
\begin{description}
\item$i)$
The transform $\mathcal T_1$ can be extended in a unique way to
a bicontinuous isomorphism 
$W^{1,p}(\R)\to \mathcal D^{1,\p}_{\rm r}(\R^n;|x|^{\alpha} dx)$.
\item$ii)$ If $\alpha>p-n$ then 
$C^1_{c,r}(\R^n)\subset \mathcal D^{1,\p}_{\rm r}(\R^n;|x|^\alpha dx)$. 
In particular, if $p<n$ then
$\mathcal D^{1,\p}_{\rm r}(\R^n;|x|^0 dx)=\mathcal D^{1,\p}_{r}(\R^n)$.
\end{description}
\end{Lemma}

\proof
Notice that
$$
\irn|x|^\alpha |\nabla(\mathcal T_1 g)|^p~dx
=\omega_n\int_{\R}\left|g'+H_\alpha~\!g\right|^\p~ds
$$
for any $g\in C^1_{c}(\R)$. Therefore $i)$ 
follows from Proposition \ref{P:equivalent1}, as $H_{\!\alpha}\neq 0$.

To prove $ii)$, take any $u\in C^1_{c,r}(\R^n)$, and let $g=\mathcal T_1^{-1} u$.
Then $g\equiv 0$ for $s<<0$ and $g(s),g'(s)=O(e^{-H_\alpha s})$ for $s\to \infty$. Since 
$H_\alpha >0$ 
then 
$g$ and $g'$ decay exponentially at infinity, and therefore $g\in W^{1,p}(\R)$. Thus 
$u\in\mathcal D^{1,p}_{\rm r}(\R^n;|x|^{\alpha})$ by $i)$, as desired.
\QED

\bigskip
\noindent
{\bf Proof of Theorem \ref{T:uniqueCKN}.} 
Using the definitions and the results in Section \ref{S:Preliminaries}, it is easy to compute
$$
S_{1,q,0}(\alpha)=
\inf_{\scriptstyle 
u\in \mathcal D^{1,\p}_{\rm r}(\R^n;|x|^\alpha dx)\atop\scriptstyle u\ne 0}
\frac{\displaystyle
\int_{\R^n}|x|^\alpha\left|\nabla u\right|^{\p} dx}
{\left(\displaystyle\int_{\R^n}|x|^{-\beta_{1,q}}|u|^q~\!dx\right)^{\p/q}}=
\omega_n^{\frac{q-\p}{q}} M_{\p,q}\!\left(-H_{\!\alpha}\right),
$$
where
$\beta_{1,q}$ is defined in (\ref{eq:beta_general}). Moreover, $u=\mathcal T_1 g$ 
solves  (\ref{eq:ELD1}) if and only if $g$ is a weak solution to
(\ref{eq:EL1}), where  $\lambda=-H_{\!\alpha}$. 
The conclusion follows by Theorem \ref{T:EF1}.
\QED

\begin{Remark}
\label{R:BW} The uniqueness result in Theorem \ref{T:uniqueCKN} has been
already stated and used in \cite{BW}, without proof.
\end{Remark}

\section{Second order inequalities}
The main results in this section concerns the infima 
$S_{2,q,0}(\alpha)$, 
$S_{2,q,1}(\alpha)$ 
and their extremals.
We follow the same scheme as in the previous section, that is, we first
prove few results about the space $W^{2,p}(\R)$;  we will turn our attention
to $\mathcal D^{2,\p}_{\rm r}(\R^n;|x|^\alpha dx)$ in subsection \ref{SS:Dk=2}.

\subsection{Equivalent norms on $W^{2,\p}(\R)$}
\label{S:W2p}

We start  
with a preliminary result. 

\begin{Lemma}
\label{L:1dim2nd}
Let $\p>1$ and $A,\gamma \in\R$ with $A^2+\gamma \ge 0$. Then
$$
I_{\p}(A,\gamma):=\inf_{\scriptstyle g\in W^{2,\p}(\R)\atop\scriptstyle g\ne 0}
\frac{\displaystyle\int_\R\left|g''-2A~\!g'-\gamma g~\!\right|^\p ds}
{\displaystyle\int_\R|g|^\p~\!ds}=|\gamma |^\p
$$
and $I_{\p}(A,\gamma)$ is not achieved. 
\end{Lemma}

\proof
Put
${\lambda}:=2+2\sqrt{A^2+\gamma }$, $b:=2+\sqrt{A^2+\gamma}-A$,
so that $(b-2)(\lambda-b)=\gamma$.
To any  $g\in C^2_c(\R)$ we associate the function
$v(r):=r^{2-b}~\!g(-\log r)$. By direct computation and using Lemma
\ref{L:1dim}
we easily infer
$$
I_{\p}(A,\gamma)=
\inf_{\scriptstyle v\in C^2_c(\R_+)\atop\scriptstyle v\ne 0}
\frac{\displaystyle
\int_0^\infty r^{pb-1}\left|v''+({\lambda}-1)r^{-1}v'\right|^\p~dr}
{\displaystyle\int_0^\infty r^{(p-2)p-1}|v|^\p~dr}=\left|(b-\lambda)(b-2)\right|^p=
|\gamma|^p~\!,
$$
and the first claim is proved. By contradiction, assume that
there exists $g\in W^{2,\p}(\R)$ such that $\displaystyle\int_\R|g|^\p~\!ds=1$
and $\displaystyle\int_\R\left|g''-2A~\!g'-\gamma g~\!\right|^\p ds=|\gamma|^p$.
Then $\gamma\neq 0$, since $g\equiv 0$ is the only function in
$W^{2,\p}(\R)$ that solves $-g''+2Ag'=0$ on $\R$. Thus we
also have that $b\neq 2$ and $b\neq \lambda$. 
Now we put $v({r})=r^{2-b}~\!g(-\log r)$ as before,
and we compute
$$
\int_0^\infty r^{(b-\lambda+1)p-1}|(r^{\lambda-1}v')'|^p~\!dr=
\int_\R\left|g''-2A~\!g'-\gamma g~\!\right|^\p~\! ds=|\gamma|^p.
$$
Then we estimate via the Hardy inequality (\ref{eq:H})
\begin{eqnarray*}
\int_0^\infty r^{(b-\lambda)p-1}|(r^{\lambda-1}v')|^p~\!dr&=&
\int_0^\infty r^{(b-1)p-1}|v'|^p~\!dr\ge
|b-2|^p\int_0^\infty r^{(b-2)p-1}|v|^p~\!dr\\
&=&
|b-2|^p\int_\R|g|^p~\!ds=|b-2|^p.
\end{eqnarray*}
This is impossible, as $\lambda-b=\frac{\gamma}{b-2}$
and since the best constant in the Hardy inequality
$$
\int_0^\infty r^{(b-\lambda+1)p-1}|\omega'|^p~\!dr\ge |\lambda-b|^\p
\int_0^\infty r^{(b-\lambda)p-1}|\omega|^p~\!dr
$$
is not achieved. The lemma is completely
proved.
\QED
Now we take an exponent $q>p$ and we study the infimum
$$
I_{\p,q}(A,\gamma):=\inf_{\scriptstyle g\in W^{2,\p}(\R)\atop\scriptstyle g\ne 0}
\frac{\displaystyle
\int_\R\left|g''-2A~\!g'-\gamma g~\!\right|^{\p} ds}
{\left(\displaystyle\int_\R|g|^q~\!ds\right)^{\p/q}}~\!.
$$
\begin{Remark}
\label{R:I=0}
If $\gamma=0$ then  $I_{\p,q}(A,0)=0$ for any $A\in\R$. For, take  
$g\in C^2_c(\R)\setminus\{0\}$ and  test $I_{\p,q}(A,0)$ with
$g_\eps(s)= g(\eps s)$, where  $\eps\to 0^+$, to get 
$$
I_{\p,q}(A,0)\le\frac{\displaystyle 
\int_\R\left| g_\eps''-2A~\!g_\eps'~\!\right|^{\p} ds}
{\left(\displaystyle\int_\R|g_\eps|^q~\!ds\right)^{\p/q}}
=\eps^{\frac{p}{q}+p-1}~\! \frac{\displaystyle 
\int_\R\left|\eps g''-2A~\!g'~\!\right|^{\p} ds}
{\left(\displaystyle\int_\R|g|^q~\!ds\right)^{\p/q}}=o(1)~\!.
$$
\end{Remark}

\begin{Proposition}
\label{P:equivalent2}
Let $\p>1$ and $A,\gamma \in\R$ with $A^2+\gamma \ge 0$. If ~$\gamma\neq 0$
then
$$
\|g\|_{A,\gamma}:=\left(\int_\R\left|g''-2A~\!g'-\gamma g~\!\right| ds\right)^{1/p}
$$
is an equivalent norm on $W^{2,\p}(\R)$. Moreover, for any $q>p$ the
infimum $I_{\p,q}(A,\gamma)$
is positive  and achieved in $W^{2,\p}(\R)$.
\end{Proposition}

\proof
Fix a small $\eps>0$ such that $|2A|\eps\le 1/2$. We recall that there exists
a constant $C_\eps>0$ such that $\|g'\|_p\le \eps\|g''\|_p+C_\eps\|g\|_p$
for any $g\in W^{2,p}(\R)$.
Using also Lemma \ref{L:1dim2nd} we find that
\begin{eqnarray*}
\|g\|_{W^{2,p}}&\le& 2
\left(1+\frac{|2A|C_\eps+|\gamma|+1}{|\gamma|}\right)
\|g''-2A~\!g'-\gamma g\|_p\\
\|g''-2A~\!g'-\gamma g\|_p&\le& 
\left(2+|2A|C_\eps+|\gamma|\right)\|g\|_{W^{2,p}}.
\end{eqnarray*}
Thus the norm
$\|\cdot\|_{A,\gamma}$ is equivalent to the standard one.

Since $W^{2,\p}(\R)\hookrightarrow L^q(\R)$ by
 Sobolev embedding, then  $I_{\p,q}(A,\gamma)>0$ by the first part of the
 proof. 
 By nowadays  standard arguments, it is easy to prove that
every bounded minimizing sequence for $I_{\p,q}(A,\gamma)$ is relatively
compact in $W^{2,p}(\R)$ up to
translations in $\R$. In particular, $I_{\p,q}(A,\gamma)$ is attained in $W^{2,\p}(\R)$.
\QED

Now we focus our attention on the inclusions $W^{2,p}(\R)\hookrightarrow W^{1,q}(\R)$,
where $q\ge p$. We start with the "linear" case $q=p$.

\begin{Lemma}
\label{L:1dim_dn}
Let $\p>1$, $A,\gamma, H \in\R$ with $A^2+\gamma \ge 0$. Let
$$
J_{\p}(A,\gamma,{H}):=\inf_{\scriptstyle g\in W^{2,\p}(\R)\atop\scriptstyle g\ne 0}
\frac{\displaystyle
\int_\R\left|g''-2A~\!g'-\gamma g~\!\right|^\p ds}
{\displaystyle\int_\R|g'+Hg|^\p~\!ds}~\!.
$$

$~~i)$ If $\gamma\neq 0$ then $J_{\p}(A,\gamma,{H})>0$.

\medskip

$~ii)$ 
$J_{\p}(A,0,0)=
\displaystyle|2A|^\p$.

\medskip

$iii)$ 
If ${H}\neq 0$ and 
${H}^2+2A{H}-\gamma=0$, then $J_{\p}(A,\gamma,{H})=\displaystyle\left|\frac{\gamma}{H}\right|^\p$.
\end{Lemma}

\proof
If $\gamma\neq 0$ then
$\|g\|_{A,\gamma}$ is an equivalent norm
on $W^{2,\p}(\R)$ by Proposition \ref{P:equivalent2}. Hence
$J_{\p}(A,\gamma,{H})>0$, since $W^{2,\p}(\R)\hookrightarrow W^{1,\p}(\R)$. 

To check $ii)$ one can reproduce the trick in the proof of Lemma \ref{L:1dim2nd}, or can 
argue as follows. First notice that $J_{\p}(A,0,0)\ge
\displaystyle|2A|^\p$ by Lemma \ref{L:1dim1st}. To prove the 
opposite inequality use a rescaling argument, as in Remark \ref{R:I=0}.
For the convenience of the reader we repeat here the proof. Take any 
$g\in C^2_c(\R)\setminus\{0\}$ and  test $J_{\p}(A,0,0)$ with
$s\mapsto g(\eps s)$, where  $\eps \to 0^+$. The conclusion is readily achieved, as
$$
J_{\p}(A,0,0)\le 
\frac{\displaystyle
\int_\R\left|\eps g''-2A~\!g'\right|^\p ds}
{\displaystyle\int_\R|g'|^\p~\!ds}=|2A|^p+o(1).
$$
It remains to check $iii)$. Notice that
$$
J_{\p}(A,\gamma,{H})=\inf_{\scriptstyle g\in W^{2,\p}(\R)\atop\scriptstyle g\ne 0}
\frac{\displaystyle
\int_\R\left|\left(g'+Hg\right)'-\frac{\gamma}{H}\left(g'+Hg\right)~\!\right|^\p ds}
{\displaystyle\int_\R|g'+Hg|^\p~\!ds} \ge 
\left|\frac{\gamma}{H}\right|^\p
$$
by Lemma \ref{L:1dim1st}. Then use rescaling as before 
to prove the opposite inequality.
\QED

In the remaining part of this section we direct our attention to "semilinear" inequalities.
For any $q>p$, $\lambda\in\R$, the infimum $M_{p,q}(\lambda)$ has been defined in (\ref{eq:M}).

\begin{Proposition}
\label{P:1dim2dn}
Let $q>\p>1$, $A,\gamma, {H}\in\R$ with $A^2+\gamma \ge 0$,
and put
$$
J_{\p,q}(A,\gamma,{H}):=\inf_{\scriptstyle g\in W^{2,\p}(\R)\atop\scriptstyle g\ne 0}
\frac{\displaystyle
\int_\R\left|g''-2A~\!g'-\gamma g~\!\right|^{\p} ds}
{\left(\displaystyle\int_\R|g'+Hg|^q~\!ds\right)^{\p/q}}.
$$

\medskip
$~~i)$ If $\gamma\neq 0$ then $J_{\p,q}(A,\gamma,{H})$ is positive and it is achieved
in $W^{2,\p}(\R)$.

$~ii)$ $J_{\p,q}(A,0,-2A)=0$ for any $A\in\R$.

$iii)$ $J_{\p,q}(A,0,0)=M_{p,q}(2A)$ and it is not achieved.
\end{Proposition}

\proof
If $\gamma\neq 0$ then
$\|g\|_{A,\gamma}$ is an equivalent norm
on $W^{2,\p}(\R)$ by Proposition \ref{P:equivalent2}. Thus $i)$ readily follows, as
$W^{2,p}(\R)\hookrightarrow W^{1,q}(\R)$. To prove that $J_{\p,q}(A,\gamma,{H})$ is attained use
standard arguments in translation-invariant problems, as for Proposition \ref{P:equivalent2}. Equality  $J_{\p,q}(A,0,-2A)=0$
can be proved via rescaling, since
$$
J_{\p,q}(A,0,-2A)=\inf_{\scriptstyle g\in W^{2,\p}(\R)\atop\scriptstyle g\ne 0}
\frac{\displaystyle
\int_\R\left|(g'-2A~\!g)'~\!\right|^{\p} ds}
{\left(\displaystyle\int_\R|g'-2Ag|^q~\!ds\right)^{\p/q}}.
$$
To check $iii)$ we first notice that
$$
J_{\p,q}(A,0,0)\ge \inf_{\scriptstyle f\in W^{1,\p}(\R)\atop\scriptstyle f\ne 0}
\frac{\displaystyle
\int_\R\left|f'-2A~\!f~\!\right|^{\p} ds}
{\left(\displaystyle\int_\R|f|^q~\!ds\right)^{\p/q}}=M_{\p,q}(2A).
$$
Next, for any  function $f\in C^1_c(\R)$, $f\neq 0$ we 
put 
$\displaystyle g(s)=\int_{-\infty}^s f(t)~\!dt$.
Then $g$ is bounded,
$g(s)\equiv 0$ for $s<<0$ and $g(s)$ is a constant for $s>>0$, so that in general
$g\notin L^p(\R)$.  Take a function $\eta\in C^2(\R)$, such that
$0\le \eta\le 1$, $\eta\equiv 1$ on $(-\infty,1)$ and $\eta\equiv 0$ on $(2,\infty)$.
We test $J_{\p,q}(A,0,0)$ with the function $g_h(s)=\eta(h^{-1}s)g(s)$, where
$h\ge 1$ is an integer. Notice that $g_h\in W^{2,p}(\R)$ since it is
smooth and it has  compact support. 
It is not difficult to show that
 $g'_h\to g'=f$ in $L^p(\R)$ and in $L^q(\R)$,  $g''_h\to g''=f'$ in $L^p(\R)$ as
 $h\to \infty$. Thus
$$
J_{\p,q}(A,0,0)\le \frac{\displaystyle
\int_\R\left|g''_h-2A~\!g'_h~\!\right|^{\p} ds}
{\left(\displaystyle\int_\R|g_h'|^q~\!ds\right)^{\p/q}}=
\frac{\displaystyle
\int_\R\left|f'-2A~\!f~\!\right|^{\p} ds}
{\left(\displaystyle\int_\R|f|^q~\!ds\right)^{\p/q}}+o(1).
$$
Thus $J_{\p,q}(A,0,0)= M_{\p,q}(2A)$, as $f$ was arbitrarily chosen.
It remains to check that $J_{\p,q}(A,0,0)$ is not attained. 
Assume that $g\in W^{2,\p}_{\rm loc}(\R)$ 
is a non constant  function such that
$g'\in L^\p(\R)\cap L^q(\R)$, $\displaystyle\int_\R|g'|^q~\!ds=1$, $g''\in L^\p(\R)$ and 
\begin{equation}
\label{eq:notachieved}
\int_\R\left|g''-2Ag'\right|^{\p} ds=J_{\p,q}(A,0,0)
= M_{p,q}(2A)~\!.
\end{equation}
Then $g'\in W^{1,\p}(\R)$ achieves $M_{\p,q}(2A)$, and hence $g'$
has constant sign (use a standard convexity argument or Theorem
\ref{T:EF1}). In particular $g$ is monotone, that implies that 
$g\notin L^\p(\R)$. Thus $g$ does not achieve
$J_{\p,q}(A,0,0)$.
\QED

\begin{Remark}
\label{R:zero}
If $A\neq 0$ then $J_{\p,q}(A,0,0)=M_{p,q}(2A)>0$ by Proposition \ref{P:equivalent1}.
In Theorem \ref{T:EF1} we proved that
the infimum $M_{\p,q}(2A)$ is achieved by a unique and positive
function
$F\in W^{1,\p}(\R)$. Therefore, any primitive
$g$ of $F$ satisfies $g', g''\in L^\p(\R)$ and (\ref{eq:notachieved}).
However, $g\notin W^{2,\p}(\R)$ since $g$ is increasing on $\R$.
\end{Remark}

In order to simplify notations we introduce the differential operators
\begin{gather*}
{\mathcal B}_+ \eta=\eta'+{H}\eta~,\quad {\mathcal B}_-=\eta'-{H}\eta~,\\
\mathcal L_+ \eta= -\eta''+ 2A\eta'+\gamma \eta~,\quad \mathcal L_-\eta=-\eta''-2A\eta'+\gamma \eta.
\end{gather*}

\begin{Remark}
\label{R:EL2}
Assume $\gamma\neq 0$. Then any minimizer $g$ for $J_{\p,q}(A,\gamma,{H})$
is, up to a Lagrange multiplier, a weak solution to the fourth order differential equation
\begin{equation}
\label{eq:EL2}
\mathcal L_-\left(|\mathcal L_+g|^{\p-2}\mathcal L_+g\right)+
{\mathcal B}_-\left(|{\mathcal B}_+ g|^{q-2}{\mathcal B}_+g\right)=0
\quad
\textrm{on $\R$.}
\end{equation}
\end{Remark}

In the next result we identify functions that coincide  
up to a change of sign and composition with translations in $\R$.

\begin{Theorem}
\label{T:uniqueJtheta}
Let $\p>1$, $q>p$, $A, {H},\gamma\in\R$ with  $\gamma, {H}\neq 0$,
$A^2+\gamma \ge 0$ and
\begin{equation}
\label{eq:BG}
{H}^2+2A{H}-\gamma=0.
\end{equation}
Then (\ref{eq:EL2})
has a unique nontrivial solution $G\in W^{2,\p}(\R)$. 
More precisely, $G$ achieves the best constant
$J_{\p,q}(A,\gamma,{H})$, 
$J_{\p,q}(A,\gamma,{H})=M_{p,q}\left(\frac{\gamma}{H}\right)$, and
\begin{equation}
\nonumber
\begin{split}
&G(s)=k~\!e^{-{H}s}\displaystyle\int_{e^{-s}}^\infty t^{\frac{\gamma}{{H}(\p-1)}-{H}-1}
\left(1+t^{\frac{\gamma(q-\p)}{{H}(\p-1)}}\right)^{\frac{\p}{\p-q}}~\!dt\quad\textit{if ${H}>0$}\\
&G(s)=k~\!
e^{-{H}s}\displaystyle\int_{0}^{e^{-s}} t^{\frac{\gamma}{{H}(\p-1)}-{H}-1}
\left(1+t^{\frac{\gamma(q-\p)}{{H}(\p-1)}}\right)^{\frac{\p}{\p-q}}~\!dt\quad\textit{if ${H}<0$,}
\end{split}
\end{equation}
where
$$
k=\left(q\left(\frac{\p}{\p-1}\right)^{\p-1}\left|\frac{\gamma}{{H}}\right|^\p\right)^{\frac{1}{q-\p}}~\!.
$$
\end{Theorem}

\proof
First of all one has to prove that $G$ is a $W^{2,p}(\R)$-solution to 
(\ref{eq:EL2}). We indicate here a way to minimize computations. We notice that
$$
G(s)=\begin{cases}
e^{-{H}s}\displaystyle\int_{-\infty}^s e^{{H}t}F(t)~\!dt&\textit{if ${H}>0$}\\
~\\
e^{-{H}s}\displaystyle\int_{s}^\infty e^{{H}t}F(t)~\!dt&\textit{if ${H}<0$,}
\end{cases}
$$
where $F\in W^{1,\p}(\R)$ is the function defined in (\ref{eq:numerare}) with 
$\lambda=\frac{\gamma}{H}$. Thus by Theorem \ref{T:EF1} we know that $F=G'+{H}G$ achieves 
the infimum $M_{\p,q}\left(\frac{\gamma}{H}\right)$ and solves 
\begin{equation}
\label{eq:copiare}
-\left(|f'\!-\!\frac{\gamma}{H}f|^{\p-2}(f'\!-\!\frac{\gamma}{H}f)\right)'\!-
\frac{\gamma}{H}|f'\!-\!\frac{\gamma}{H}f|^{\p-2}(f'\!-\!\frac{\gamma}{H}f)=
|f|^{q-2}f.
\end{equation}
Since $G$ decays exponentially at $\pm\infty$, then
clearly $G\in L^{\p}(\R)$. Hence $G\in W^{2,\p}(\R)$,
as $G'=F-{H}G\in L^p(\R)$.
Now we use (\ref{eq:BG}) to get 
$$F'-\frac{\gamma}{H}F=G''-2AG'-\gamma G= -\mathcal L_+G.$$ 
Hence, we have showed
that $G$ solves
\begin{eqnarray}
\nonumber
\left(|\mathcal L_+G|^{\p-2}(\mathcal L_+G)\right)'+
\frac{\gamma}{H}|\mathcal L_+G|^{\p-2}(\mathcal L_+G)&=&
\label{eq:15}
|G'+{H}G|^{q-2}(G'+{H}G)\\
&=&|{\mathcal B}_+G|^{q-2}{\mathcal B}_+G~\!.
\end{eqnarray}
Finally, we apply the operator $-{\mathcal B}_-$ to both sides of (\ref{eq:15})
and we use again (\ref{eq:BG})
to get that $G$ is a solution to (\ref{eq:EL2}).

Now, assume that $g\in W^{2,\p}(\R)\setminus\{0\}$ solves (\ref{eq:EL2}), and put
$$
f:=g'+Hg={\mathcal B}_+g\in W^{1,p}(\R).
$$
We have to show that $f$ solves (\ref{eq:copiare}).
By (\ref{eq:BG}) we have that
$$
\f:=\left|f'-\frac{\gamma}{H}f\right|^{\p-2}\left(f'-\frac{\gamma}{H}f\right)=
-|\mathcal L_+g|^{p-2}\mathcal L_+g\in L^{p'}(\R)
$$
with pointwise a.e. equalities. Since $g$ solves (\ref{eq:EL2}) then
$\f$ is a distributional solution to
\begin{equation}
\label{eq:fiB}
\mathcal L_-\f={\mathcal B}_-\left(|f|^{q-2}f\right),
\end{equation}
that is,
$$
\int_\R\f(\mathcal L_+\eta)~\!ds=-\int_\R |f|^{q-2}f(\mathcal B_+\eta)~\!ds
\quad\textrm{for any $\eta\in W^{2,p}(\R)$}
$$
 (use a density argument).
Now take
any $w\in C^\infty_c(\R)$ and  put
$$
\eta(s)=
\begin{cases}
~~\!e^{-Hs}\displaystyle\int_{-\infty}^s e^{Ht}w(t)~\!dt&\textrm{if $H>0$}\\
~&\\
-~\!e^{-Hs}\displaystyle\int_{s}^{+\infty} e^{Ht}w(t)~\!dt&\textrm{if $H<0$.}
\end{cases}
$$
Notice that $\eta\in W^{2,p}(\R)$ since $H\neq 0$ and
$\eta'+H\eta=w$. Moreover, it holds that 
$$
\int_\R\left|\mathcal B_+\eta\right|^p~\!ds\le c\|\eta\|^p_{W^{1,p}}\le c\int_\R|w|^p~\!ds~,\quad
\mathcal L_+\eta=-w'+\frac{\gamma}{H}w
$$
by (\ref{eq:BG}), and 
$|f|^{q-1}\in L^{\p'}(\R)$ as $g\in W^{2,\p}(\R)\hookrightarrow W^{1,\tau}(\R)$
for any $\tau\ge \p$. Thus
\begin{eqnarray*}
\left|\int_\R\f w'~\!ds\right|&\le&
\left|\int_\R\f~\!(-\mathcal L_+\eta+\frac{\gamma}{H}~\!
w)\right|~\!ds\\
&\le&
\left|\int_\R |f|^{q-2}f\left(\mathcal B_+\eta\right)\right|+
c\left(\int_\R|w|^p~\!ds\right)^{1/p}
\le c \left(\int_\R|w|^p~\!ds\right)^{1/p}.
\end{eqnarray*}
Hence $\f\in W^{1,\p'}(\R)$
and $\f$ is a weak solution to (\ref{eq:fiB}). On the other hand, in the dual
$W^{-1,\p}(\R)$ we can compute
\begin{eqnarray*}
\mathcal L_-\f&=&
-\f''-2A\f'+\gamma\f\\
&=&-\left(\f'+\frac{\gamma}{H}\f\right)'+
{H}\left(\f'+\frac{\gamma}{H}\f\right)=-{\mathcal B}_-\left(\f'+\frac{\gamma}{H}\f\right),
\end{eqnarray*}
thanks to (\ref{eq:BG}).
Thus
we have shown that
$$
-{\mathcal B}_-\left(\f'+\frac{\gamma}{H}\f\right)=
{\mathcal B}_-\left(|f|^{q-2}f\right)~\!.
$$
The operator 
${\mathcal B}_-:W^{1,\p}(\R)\to W^{-1,\p}(\R)$ is invertible, and therefore
it holds that $-\f'-\frac{\gamma}{H}\f=|f|^{q-2}f$, that is, 
$f$ solves (\ref{eq:copiare}).
By Theorem \ref{T:EF1} we can assume that $f$ coincides with the function $F=G'+{H}G$.
Thus $g'+Hg=G'+{H}G$, that implies $g=G$, as $g,G\in W^{1,p}(\R)$.
The theorem is completely proved.
\QED

\subsection{The space $\mathcal D^{2,\p}_{\rm r}(\R^n;|x|^\alpha dx)$}
\label{SS:Dk=2}
In order to simplify
notation we put
$$
{H}_2={H}_{\alpha,2}= \frac{n+\alpha}{\p}-2~,
\quad 
\gamma_{2}=\gamma_{\alpha,2}=
\left(\frac{n+\alpha}{p}-2\right)\left(n-\frac{n+\alpha}{p}\right),
$$
compare with (\ref{eq:B}) and (\ref{eq:BGamma}). We need also the constant
$$
A_{2}= \frac{n-2}{2}-{H}_2~\!.
 $$
 Notice that
\begin{equation}
\label{eq:BGR}
A_{2}^2+\gamma_{2}=\left(\frac{n-2}{2}\right)^2\ge 0~,\quad
{H}_{2}^2+2A_{2}{H}_{2}-\gamma_{2}=0.
\end{equation}

Assume that $\gamma_2\neq 0$, 
that is, $\alpha\notin\{2p-n,np-n\}$. 
Then $\mathcal D^{2,\p}_{\rm r}(\R^n;|x|^\alpha dx)$ 
is a Banach space with norm
$$
\|u\|_{2,p}^p=\irn|x|^\alpha|\Delta u|^p~dx~.
$$
Moreover, $\mathcal D^{2,\p}_{\rm r}(\R^n;|x|^\alpha dx)$ is continuously embedded into $L^p(\R^n;|x|^{\alpha-2p}dx)$ and into
$\mathcal D^{1,\p}_{\rm r}(\R^n;|x|^{\alpha-p} dx)$
by (\ref{eq:R}), (\ref{eq:AS1}).

For $g\in C^2_c(\R)$ we put $(\mathcal T_2g)(x)=|x|^{-{H}_{2}}g(-\log|x|)$.
Then $\mathcal T_2$ is a linear, invertible transform 
$C^2_c(\R)\to C^2_{c,r}(\Rbuco)$.

Now we prove the second-order version of Lemma \ref{L:space_explicit1}.
\begin{Lemma}
\label{L:space_explicit2}
Assume $\alpha\notin\{2p-n,np-n\}$.
\begin{description}
\item$i)$
The transform $\mathcal T_2$ can be extended in a unique way to
a bicontinuous isomorphism 
$W^{2,p}(\R)\to \mathcal D^{2,\p}_{\rm r}(\R^n;|x|^{\alpha} dx)$.
\item$ii)$ If $\alpha>2p-n$ then 
$C^2_{c,r}(\R^n)\subset \mathcal D^{2,\p}_{\rm r}(\R^n;|x|^\alpha dx)$. 
In particular, if $n>2p$ then
$\mathcal D^{2,\p}_{\rm r}(\R^n;|x|^0 dx)=\mathcal D^{2,\p}_{r}(\R^n)$.
\end{description}
\end{Lemma}

\proof
By direct computation one can check that 
$$
\int_{\R^n}|x|^{\alpha}|\Delta (\mathcal T_2 g)|^\p dx=\omega_n\int_{\R}
\left|g''-2A_{2}~\!g'-\gamma_{2}g\right|^\p~ds
$$
for any $g\in C^2_{c,r}(\R)$. Therefore Proposition \ref{P:equivalent2}
immediately implies
$i)$. To prove $ii)$ fix $u\in C^2_{c,r}(\R^n)$,  put $g=\mathcal T_2^{-1} u$
and then argue as in Lemma \ref{L:space_explicit1}.
\QED

Now we fix an exponent $q>\p$ and we use Lemma \ref{L:space_explicit2}
together with the results in Section \ref{S:W2p} to study the best constants
$S_{2,q,0}(\alpha)$ and $S_{2,q,1}(\alpha)$. 
Notice that 
\begin{gather*}
\int_{\R^n}|x|^{-\beta_{2,q}}|\mathcal T_2 g|^q~\!dx= \omega_n
\int_{\R}|g|^q ds\\
\int_{\R^n}|x|^{-\beta_{1,q}}|\nabla(\mathcal T_2 g)|^q~\!dx=\omega_n\int_{\R}
\left|g'+{H}_2~\!g\right|^q~ds
\end{gather*}
for any 
$g\in W^{2,p}(\R)$ 
where, accordingly with (\ref{eq:beta_general}),
$$
\beta_{2,q}=n-q\frac{n-2\p+\alpha}{\p}~,\quad
\beta_{1,q}=n-q\frac{n-\p+\alpha}{\p}~\!.
$$

First we use Lemma \ref{L:space_explicit2} and the above computations
to observe  that the minimization problems
\begin{gather*}
S_{2,q,0}(\alpha)=
\inf_{\scriptstyle u\in\mathcal D^{2,\p}_{\rm r}(\R^n;|x|^\alpha dx)\atop\scriptstyle u\ne 0}
\frac{\displaystyle
\int_{\R^n}|x|^{\alpha}|\Delta u|^\p dx}
{\left(\displaystyle\int_{\R^n}|x|^{-\beta_{2,q}}| u|^q dx\right)^{\p/q}}\\
I_{\p,q}(A_{2},\gamma_{2})=
\inf_{\scriptstyle g\in W^{2,\p}(\R)\atop\scriptstyle g\ne 0}
\frac{\displaystyle
\int_\R\left|g''-2A_2~\!g'-\gamma_2 g~\!\right|^{\p} ds}
{\left(\displaystyle\int_\R|g|^q~\!ds\right)^{\p/q}}
\end{gather*}
are equivalent. Thus  Proposition \ref{P:equivalent2}
immediately implies the next existence result, which is
indeed  included in the more general Theorem \ref{T:extremalsk}.

\begin{Theorem}
\label{T:RSachievedI}
If $\alpha\notin\{2p-n,np-n\}$ and $q>p$,  then
$S_{2,q,0}(\alpha)$ is positive and attained 
 in $\mathcal D^{2,\p}_{\rm r}(\R^n;|x|^\alpha dx)$.
\end{Theorem}
We notice that any minimizer for $S_{2,q,0}(\alpha)$ satisfies the Euler-Lagrange equation
$$
\Delta\left(|x|^\alpha|\Delta u|^{p-2}\Delta u\right)=|x|^{-\beta_{2,q}}| u|^{q-2} u~\!,
$$
which is equivalent to the well known and largely studied H\'enon-Lane-Emden system
$$
\begin{cases}
-\Delta u=|x|^{a}|v|^{\p'-2}v\\
-\Delta v=|x|^{b}|u|^{q-2}u,
\end{cases}
$$
where $a=-\frac{\alpha}{\p-1}$ and $b=-\beta_{2,q}$.
In particular Theorem \ref{T:RSachievedI} provides the existence of solutions
to the the above system, if 
 $a,b\neq -n$ and $p',q$ lie on the "{critical hyperbola}"
 $$
 \frac{a+n}{p'}+\frac{b+n}{q}=n-2~\!.
 $$
 We quote  \cite{MS} for details and additional results.
 
We are in position to prove one of the main results in the introduction.

\bigskip
\noindent
{\bf Proof of Theorem \ref{T:RSachievedJ}.}
For any $f\in W^{2,p}(\R)$ put
$$
\mathcal L f=-f''+2A_2f+\gamma_2 f~,\quad \mathcal B f=f'+H_2 f.
$$
Direct computations lead to 
$$
-\Delta\left(\mathcal T_2 f\right)=\mathcal T_0\left(\mathcal L f\right)~,\quad
\left(\mathcal T_2 f\right)'=-\mathcal T_1\left(\mathcal B f\right),
$$
see (\ref{eq:EF}) for the definition of $\mathcal T_h$.
In particular, for any pair $U=\mathcal T_2 G$, $\f=\mathcal T_2 f$
we have that
\begin{gather*}
\irn|x|^\alpha|\Delta U|^{p-2}\Delta U\Delta\f~\!dx=
\omega_n \int_\R
|\mathcal L G|^{p-2}(\mathcal L  G)(\mathcal L f)~\!ds\\
\irn|x|^{-\beta_{1,q}}|\nabla U|^{p-2}\nabla U\cdot\nabla\f~\!dx=
\omega_n \int_\R
|\mathcal B G|^{q-2}(\mathcal B  G)(\mathcal B f)~\!ds.
\end{gather*}
To conclude the proof use Theorem \ref{T:uniqueJtheta} with $A=A_2$ and
$$
\frac{\gamma}{H}=\frac{\gamma_2}{H_2}=\frac{np-n+\alpha}{p}.
$$
\QED

\begin{Remark}
\label{R:degenerate1}
Here we take $\alpha=np-n$ and we define
\begin{gather*}
S_{2,q,0}(np-n):=
\inf_{\scriptstyle u\in C^2_{c,r}(\R^n\!\setminus\!\{0\})\atop\scriptstyle u\ne 0}
\frac{\displaystyle
\int_{\R^n}|x|^{np-n}|\Delta u|^\p dx}
{\left(\displaystyle\int_{\R^n}|x|^{-n+q(n-2)}| u|^q dx\right)^{\p/q}}\\
S_{2,q,1}(np-n):=
\inf_{\scriptstyle u\in C^2_{c,r}(\R^n\!\setminus\!\{0\})\atop\scriptstyle u\ne 0}
\frac{\displaystyle
\int_{\R^n}|x|^{np-n}|\Delta u|^\p dx}
{\left(\displaystyle\int_{\R^n}|x|^{-n+q(n-1)}| \nabla u|^q dx\right)^{\p/q}}~\!.
\end{gather*}
Thanks to the Emden-Fowler transform
$(\mathcal T_2 g)(x)=|x|^{2-n}g(-\log|x|)$, and using 
Remark \ref{R:I=0} and $ii)$ in Proposition \ref{P:1dim2dn},
it is easy to check that
\begin{gather*}
S_{2,q,0}(np-n)=\omega_n^{\frac{q-p}{p}}I_{p,q}\left(-\frac{n-2}{2},0\right)=0\\
S_{2,q,1}(np-n)=\omega_n^{\frac{q-p}{p}}J_{p,q}\left(-\frac{n-2}{2},0,n-2\right)=0.
\end{gather*}
\end{Remark}

\begin{Remark}
\label{R:degenerate2}
Here we take $n\ge 3$ and $\alpha=2p-n$ (notice that for $n=2$ the two 
"degenerate cases" $\alpha=np-n$ and $\alpha=2p-n$ coincide).
Now we observe that $H_2=0$ and we use $(\mathcal T_2 g)(x)=g(-\log|x|)$ to get
$$
S_{2,q,0}(2p-n):=
\inf_{\scriptstyle u\in C^2_{c,r}(\R^n\!\setminus\!\{0\})\atop\scriptstyle u\ne 0}
\frac{\displaystyle
\int_{\R^n}|x|^{2p-n}|\Delta u|^\p dx}
{\left(\displaystyle\int_{\R^n}|x|^{-n}| u|^q dx\right)^{\p/q}}=
\omega_n^{\frac{q-p}{p}}I_{p,q}\left(\frac{n-2}{2},0\right)=0
$$ 
by  Remark \ref{R:I=0}. Next we define
$$
S_{2,q,1}(2p-n):=
\inf_{\scriptstyle u\in C^2_{c,r}(\R^n\!\setminus\!\{0\})\atop\scriptstyle u\ne 0}
\frac{\displaystyle
\int_{\R^n}|x|^{2p-n}|\Delta u|^\p dx}
{\left(\displaystyle\int_{\R^n}|x|^{-n+q}| \nabla u|^q dx\right)^{\p/q}}~\!.
$$
We have that
$$
S_{2,q,1}(2p-n)=\omega_n^{\frac{q-p}{p}}J_{p,q}\left(\frac{n-2}{2},0,0\right)=M_{p,q}(n-2)
$$
by $iii)$ in Proposition  \ref{P:1dim2dn}. Thus the value of
$S_{2,q,1}(2p-n)$
is explicitly known, thanks to Theorem \ref{T:EF1}. In particular, 
$S_{2,q,1}(2p-n)>0$, but 
no function in $L^p(\R^n;|x|^{-n}dx)$ achieves $S_{2,q,1}(2p-n)$,
use again Proposition \ref{P:1dim2dn}. 
See also Remark \ref{R:critico} below for the case $\alpha=2p-n=0$, $q=p^*$.
\end{Remark}

\begin{Remark}
\label{R:L1loc}
Assume  $\alpha>2p-n$, $\alpha\neq n(p-1)$. Then the weight
$|x|^{\alpha-2p}$ is locally integrable and the space $C^2_{c,r}(\R^n)$
is dense in $\mathcal D^{2,\p}_{\rm r}(\R^n;|x|^\alpha dx)$ by $ii)$ in Lemma
\ref{L:space_explicit2}. In particular, 
inequalities (\ref{eq:R}) and (\ref{eq:AS1}) hold in 
$C^2_{c,r}(\R^n)$. We also have that
\begin{gather*}
\int_{\R^n}|x|^{\alpha}|\Delta u|^\p dx\ge S_{2,q,0}(\alpha)
\left(\displaystyle\int_{\R^n}|x|^{-n+q\frac{n-2p+\alpha}{p}}| u|^q dx\right)^{\p/q}
\quad\forall u\in C^2_{c,r}(\R^n)\\
\int_{\R^n}|x|^{\alpha}|\Delta u|^\p dx\ge S_{2,q,1}(\alpha)
\left(\displaystyle\int_{\R^n}|x|^{-n+q\frac{n-p+\alpha}{p}}|\nabla u|^q dx\right)^{\p/q}
\quad\forall u\in C^2_{c,r}(\R^n).
\end{gather*}
\end{Remark}

\begin{Remark}
\label{R:critico}
Take $\alpha=0$. If  $n=2p$ then $p^*=n$, the  constant
$S_{2,n,1}(0)$ in Remark \ref{R:degenerate2} 
is positive and  any radially symmetric function $\widetilde U$  satisfying
$$
\widetilde U'({r})=kr(1+r^n)^{-1}~,\quad k=2p\left(2(p-1)\right)^{1/p}=n(n-2)^{2/n}
$$ 
is indeed a solution to the Euler-Lagrange equation (\ref{eq:ELSobolev}). For,
it is convenient to notice that 
$\Delta\left(|\Delta \widetilde U|^{p-2}\Delta \widetilde U\right)+\Delta_{p^*} \widetilde U=
\div\left(\Phi|x|^{-1}x\right)$, where
$$
\Phi:=\left(|\Delta \widetilde U|^{p-2}\Delta \widetilde U\right)'+
|\widetilde U'|^{p^*-2}\widetilde U'.
$$
We compute 
$\Delta \widetilde U=nk(1+r^{n})^{-2}=n^2(n-2)^{2/n}(1+r^{n})^{-2}$, so that
$$
\left(|\Delta \widetilde U|^{p-2}\Delta \widetilde U\right)'=
-n^{n-1}(n-2)^{\frac{2(n-1)}{n}}~\!r^{n-1}(1+r^n)^{1-n}
=-|\widetilde U'|^{p^*-2}\widetilde U'.
$$
Thus $\Phi\equiv 0$ and therefore $\widetilde U$ solves  (\ref{eq:ELSobolev}). We also point out that
$$
\irn|\Delta \widetilde U|^pdx<\infty~,\quad\irn|x|^{-p}|\nabla \widetilde U|^pdx<\infty~,
\quad \irn|x|^{-2p}| \widetilde U|^pdx=\infty,
$$
as $\widetilde U$ is radially increasing and the weight
$|x|^{-2p}=|x|^{-n}$ is not integrable at $0$, nor at infinity.
For instance, for any constant $c\in\R$ the function
$ \widetilde U(x)=2\sqrt 2 \arctan |x|^{2}+c$ solves 
$$
\Delta^{\!2} U+\Delta_{4} U=0\quad\textit{on $\R^4$}.
$$
\end{Remark}

\section{Higher order inequalities}
In this last section we study the case $k>2$ and we prove Theorem \ref{T:extremalsk}.
We start by investigating the standard space $W^{k,p}(\R)$.

\subsection{Equivalent norms on $W^{k,\p}(\R)$}
\label{S:higher}
Let $m\ge 1$ be an integer and let
$$
\vec A=(A_1,...,A_m)\in\R^m~,\quad \vec \gamma=(\gamma_1,...,\gamma_m)\in\R^m
$$
be given $m$-vectors. We define the $m+1$ differential operators
$$
L_h g=g''-2A_{h}~\!g'-\gamma_{h}~\!g~,\quad \mathbb L_{\vec A,\vec \gamma}=L_1\circ...\circ L_m~\!,
$$
so that $\mathbb L_{\vec A,\vec \gamma}$ has order $2m$.

We distinguish the "even case" $k=2m$ from the "odd" one, when $k=2m+1$.

\begin{Proposition}  
\label{P:equivalent2m}
Assume that $A_{h}^2+\gamma_{h} \ge 0$ for any $h=1,...,m$.
Then
$$
I_{\p}(\vec A,\vec \gamma):=\inf_{\scriptstyle g\in W^{2m,\p}(\R)\atop\scriptstyle g\ne 0}
\frac{\displaystyle
\int_\R\left|\mathbb L_{\vec A,\vec \gamma} ~\!g~\!\right|^\p ds}
{\displaystyle\int_\R|g|^\p~\!ds}=\prod_{h=1}^m|\gamma_{h} |^{\p}.
$$
Moreover, $\|g\|_{\vec A,\vec \gamma}:=
\|\mathbb L_{\vec A,\vec \gamma} ~\!g\|_\p$ is an equivalent norm on $W^{2m,p}(\R)$
provided that $\gamma_{h}\neq 0$ for any $h=1,...,m$. 
\end{Proposition}

\proof
To check that
$I_{\p}(\vec A,\vec \gamma)\le \displaystyle\prod_{h=1}^m|\gamma_{h} |^{\p}$
we use  a rescaling
 argument, as in Lemma \ref{L:1dim_dn}. Take any 
$g\in C^k_c(\R)\setminus\{0\}$ and  test $I_{\p}(\vec A,\vec \gamma)$ with
$g_\eps(s):=g(\eps s)$, where  $\eps \to 0^+$.  
By direct computation one gets that
$$
\int_\R\left|\mathbb L_{\vec A,\vec \gamma} ~\!g_\eps~\!\right|^\p ds=
\eps^{-1}\!\left(~\!\prod_{h=1}^m|\gamma_{h} |^{\p}\right)\int_\R|g|^p~\!ds+o(\eps^{-1})~,\quad
\int_\R|g_\eps|^\p~\!ds=\eps^{-1} \!\int_\R|g|^\p~\!ds.
$$
The conclusion readily follows.
 The opposite inequality can be proved by induction, starting with the case $k=2$
 that has already been discussed in Section \ref{S:W2p}.
 \QED
 
 The same arguments lead to the next result. We omit the details of the proof.

\begin{Proposition}  
\label{P:equivalent2m+1}
Assume that $A_{h}^2+\gamma_{h} \ge 0$ for any $h=1,...,m$ and let $\lambda\in\R$. 
Then 
$$
M_{\p}(\vec A,\vec \gamma;\lambda):=\inf_{\scriptstyle g\in W^{2m+1,\p}(\R)\atop\scriptstyle g\ne 0}
\frac{\displaystyle
\int_\R\left|\mathbb L_{\vec A,\vec \gamma} ~\!g'-\lambda\mathbb L_{\vec A,\vec \gamma} ~\!g~\!\right|^\p ds}
{\displaystyle\int_\R|g|^\p~\!ds}=|\lambda|^\p\prod_{h=1}^m|\gamma_{h} |^{\p}.
$$
Moreover, $\|g\|_{\vec A,\vec \gamma,\lambda}:=
\|\mathbb L_{\vec A,\vec \gamma} ~\!g'-\lambda\mathbb L_{\vec A,\vec \gamma} ~\!g\|_\p$ is an equivalent norm on $W^{2m+1,p}(\R)$ provided that $\lambda\neq 0$ and 
$\gamma_{h}\neq 0$ for any $h=1,...,m$.
\end{Proposition}

\begin{Remark}
One can get more inequalities by taking advantage of the embeddings
$W^{k,p}(\R)\hookrightarrow W^{j,q}(\R)$ for $h=1,...,k-1$ and $q\ge p$.
\end{Remark}

\subsection{The space $\mathcal D^{k,\p}_{\rm r}(\R^n;|x|^\alpha dx)$}
\label{A:poli}
To simplify notation we put
$$
{H}_h={H}_{\alpha,h}= \frac{n+\alpha}{\p}-h~,
\quad 
\gamma_{h}=\gamma_{\alpha,h}=
{H}_h\left(n-2-{H}_h\right)
$$
for any integer $h\ge 1$, compare with (\ref{eq:B}) and (\ref{eq:BGamma}).
We introduce also the constants
$$
A_{h}= \frac{n-2}{2}-{H}_h.
$$
Notice that 
$$
A_{h}^2+\gamma_h=\left(\frac{n-2}{2}\right)^2\ge 0~,\quad
{H}_h^2+2A_{h}{H}_h-\gamma_h=0.
$$

In this section we will always assume that (\ref{eq:nondeg_gen}) is satisfied. In particular, we have that 
$\mathcal D^{k,\p}_{\rm r}(\R^n;|x|^\alpha dx)$ is a well defined Banach space with norm
$$
\|u\|_{k,\alpha}^p=\irn|x|^\alpha|\nabla^k u|^p~dx~\!,
$$
and  it is continuously embedded into $L^p(\R^n;|x|^{\alpha-kp}dx)$
by Theorem \ref{T:theorem_gen}. 

Accordingly with   (\ref{eq:EF}), we put $(\mathcal T_kg)(x)=|x|^{-{H}_{k}}g(-\log|x|)$ for  $g\in C^k_c(\R)$.

\begin{Lemma}  
\label{L:RRSk}
The
transform
$\mathcal T_{k}$  can be extended in a unique way to
a bicontinuous isomorphism 
$W^{k,p}(\R)\to \mathcal D^{k,\p}_{\rm r}(\R^n;|x|^{\alpha} dx)$. 
\end{Lemma}

\proof
Recall that $(\mathcal T_h g)(x)=|x|^{-\frac{n+\alpha}{p}+h}g(-\log |x|)$
for any $h\ge 0$ and
notice that
\begin{equation}
\label{eq:induction}
\Delta(\mathcal T_{h}g)=\mathcal T_{h-2}\left(g''-2A_{h}~\!g'-\gamma_h~\! g\right)
\end{equation}
for any $h\ge 2$ (use an induction argument).

We distinguish the case $k=2m$ from the case of odd order operators.

\bigskip
\noindent
{\bf Even poliharmonic operators.}
Assume that $k=2m$ is an even integer and that
$\gamma_{2h}\neq 0$ for any $h=1,...,m$.
We adopt the notation in Section \ref{S:higher} with
\begin{gather*}
\vec A=(A_2,..,A_{2m})~,\quad \vec\gamma=(\gamma_2,...,\gamma_{2m}),\\
L_hg = g''-2A_{2h}~\!g'-\gamma_{2h}~\!g~,\quad \mathbb L_{\vec A,\vec \gamma}=L_1\circ...\circ L_m.
\end{gather*}
Using (\ref{eq:induction}) it is easy to prove by induction that
$$
\Delta^m ( \mathcal T_{2m}g) = \mathcal T_0\left(\mathbb L_{\vec A,\vec \gamma} ~\!g\right)
$$
for any $g\in C^k_c(\R)$.
Therefore, for $u= \mathcal T_{2m}(g)$ the
following equality holds
$$
\|u\|^\p_{2m,\alpha}=\int_{\R^n}|x|^{\alpha}|\Delta^m u|^\p dx=\omega_n\int_{\R}
\left|\mathbb L_{\vec A,\vec \gamma}g\right|^\p~ds.
$$
The conclusion in the even case follows by Proposition
\ref{P:equivalent2m}.

\bigskip
\noindent
{\bf Odd poliharmonic operators.}
When $k=2m+1$ is odd we assume that
$\gamma_{2h+1}\neq 0$ for any $h=1,...,m$ and that the Hardy constant is positive.
Now we put 
\begin{gather*}
\vec A=(A_3,..,A_{2m+1})~,\quad \vec\gamma=(\gamma_3,...,\gamma_{2m+1}),\\
L_hg = g''-2A_{2h+1}~\!g'-\gamma_{2h+1}~\!g~,\quad \mathbb L_{\vec A,\vec \gamma}=L_1\circ...\circ L_m.
\end{gather*}
Using Ä(\ref{eq:induction}) one can prove by induction that
$$
\Delta^m ( \mathcal T_{2m+1}g)=\mathcal T_{1}(\mathbb L_{\vec A,\vec \gamma} ~\!g)
$$
for any $g\in C^k_c(\R)$.
Therefore, for $u= \mathcal T_{2m+1}g$ it is easy to compute
$$
\|u\|^\p_{2m+1,\alpha}=\int_{\R^n}|x|^{\alpha}|\nabla(\Delta^m u)|^\p dx=
\omega_n\int_{\R}
\left|\mathbb L_{\vec A,\vec \gamma}~\!g'+H_{\!\alpha}\mathbb L_{\vec A,\vec \gamma}~\!g \right|^\p~ds.
$$
The conclusion follows by Proposition
\ref{P:equivalent2m+1}.
\QED

\begin{Remark}
\label{R:alternative}
The computations in the proof of Lemma \ref{L:RRSk} together with Propositions 
\ref{P:equivalent2m}, \ref{P:equivalent2m+1} provide an alternative proof
of Theorem \ref{T:theorem_gen}.
\end{Remark}

\begin{Remark}
\label{R:densityk}
Assume $\alpha>kp-n$.
One can adapt the arguments already used in the lowest order case $k=1$
to show that 
$C^k_{c,r}(\R^n)\subset \mathcal D^{k,\p}_{\rm r}(\R^n;|x|^\alpha dx)$.
\end{Remark}

\begin{Remark}
\label{R:regularity}
Thanks to Lemma \ref{L:RRSk} we can identify the spaces
$\mathcal D^{k,\p}_{\rm r}(\R^n;|x|^\alpha dx)$ and
$W^{k,p}(\R)$ trough $\mathcal T_k$. In particular,
every function $u\in \mathcal D^{k,\p}_{\rm r}(\R^n;|x|^\alpha dx)$
has continuos derivatives up to the order $k-1$, and partial derivatives of order $k$ 
(in the classical sense) exist for almost every $|x|>0$. 
\end{Remark}

Now we focus our attention on Sobolev type embeddings. For 
$j\in\{0,...,k\}$ and $q>p$ define $\beta_{k-j,q}$ accordingly with
(\ref{eq:beta_general}).

\begin{Remark}
\label{R:Dj}
Fix an index $j=1,..., k-1$, and notice that
$$
\widetilde {H}_h:= \frac{n-\beta_{k-j,q}}{q}-h=\frac{n+\alpha}{p}-(k-j+h)={H}_{k-j+h},
$$
so that $\widetilde {H}_j={H}_k$. Moreover, $\widetilde\gamma_h:=\widetilde {H}_h(n-2+\widetilde {H}_h)$ satisfies
(\ref{eq:nondeg_gen}) with $\alpha,k,p$ replaced by $\beta_{k-j,q},j$ and $q$, respectively.
Thus $\mathcal T_k$ can be used to identify $W^{j,q}(\R)$ with the Banach space
$\mathcal D^{j,q}_{\rm r}(\R^n;|x|^{-\beta_{k-j,p}} dx)$ by Lemma 
\ref{L:RRSk}.
\end{Remark}

We are in position to prove Theorem \ref{T:extremalsk}.

\bigskip
\noindent
{\bf Proof of Theorem \ref{T:extremalsk}.}
By Lemma \ref{L:RRSk}, Sobolev embedding theorem and Remark \ref{R:Dj}, we
have the following chain of continuous arrows:
$$
\mathcal D^{k,\p}_{\rm r}(\R^n;|x|^\alpha dx)\stackrel{\mathcal T_k^{-1}}
{\longrightarrow} W^{k,p}(\R) 
{\hookrightarrow} W^{j,q}(\R) \stackrel{\mathcal T_k}
{\longrightarrow} \mathcal D^{j,q}_{\rm r}(\R^n;|x|^{-\beta_{k-j,q}} dx).
$$
Thus  $\mathcal D^{k,\p}_{\rm r}(\R^n;|x|^\alpha dx)\hookrightarrow
\mathcal D^{j,q}_{\rm r}(\R^n;|x|^{-\beta_{k-j,q}} dx)$
with a continuous embedding,
and hence $S_{k,q,j}(\alpha)$ is positive.

To prove that $S_{k,q,j}(\alpha)$ is achieved one can study an equivalent 
minimization problem for functions in $W^{k,p}(\R)$, as we did in Section \ref{SS:Dk=2}
for the case $k=2$. Here we adopt a direct approach. 
The strategy is essentially the same of \cite{Mu0}, and it was 
originally inspired to the author by the famous paper
\cite{SaUh} by Sacks and Uhlenbeck. 

We introduce the energies $e_{k,p}, e_{j,q}$ and  $E:\mathcal D^{k,\p}_{\rm r}(\R^n;|x|^\alpha dx)\to \R$
by putting
\begin{gather*}
e_{k,p}(u)={\frac{1}{p}}\irn|x|^{\alpha}|\nabla^k u|^p~\!dx~,
\quad e_{j,q}(u)={\frac{1}{q}}\irn|x|^{-\beta}|\nabla^j u|^q~\!dx,\\
E(u)=e_{k,p}(u)- e_{j,q}(u),
\end{gather*}
where $\beta=\beta_{k-j,q}$.
The functional $E$ is of class $C^1$ and 
minimizers for $\widetilde S:=S_{k,q,j}(\alpha)$ give rise to critical points
of $E$.
Using Ekeland's variational principle
we can select a minimizing sequence 
$u_h$ such that
\begin{equation}
\label{eq:energy}
E'(u_h)\cdot v\to 0
\quad\textrm{uniformly for $v$ in bounded subsets of 
$\mathcal D^{k,\p}_{\rm r}(\R^n;|x|^\alpha dx)$.}
\end{equation}
In particular we have  
 $o(\|u_h\|_{k,\alpha})=E'(u_h)\cdot u_h=\|u_h\|_{k,\alpha}^p+O(\|u_h\|_{k,\alpha}^q)$.
Thus $u_h$ is bounded and we may suppose, without loss of generality, that
$$
\irn|x|^{-\beta}|\nabla^j u_h|^q~\!dx=\irn|x|^{\alpha}|\nabla^k u_h|^p~\!dx+o(1)=
\widetilde S^{\frac{q}{q-p}}+o(1).
$$
Since the ratio in (\ref{eq:theorem_semilinear}) is invariant with respect to dilations,
we can assume that
\begin{equation}
\label{eq:rescaling}
\int_{\{|x|<2\}}|x|^{-\beta}|\nabla^j u_h|^q~\!dx=\left(\frac{1}{2}~\!\widetilde S\right)^{\frac{q}{q-p}}.
\end{equation}
We have to prove that, up to a subsequence, $u_h$ converges weakly to some
nontrivial limit $u$. Then, a standard convexity argument shows that 
$u$ achieves $\tilde S$. 
Assume by contradiction that $u_h\weak 0$ weakly in
$\mathcal D^{k,\p}_{\rm r}(\R^n;|x|^\alpha dx)$.  By Lemma \ref{L:Rellich} we have that $|\nabla^j u_h|\to 0$ in $L^q_{\rm loc}(\Rbuco)$
and therefore    
\begin{equation}
\label{eq:contradiction}
\int_{\{|x|<1\}}|x|^{-\beta}|\nabla^j u_h|^q~\!dx=
\left(\frac{1}{2}~\!\widetilde S\right)^{\frac{q}{q-p}}+o(1)
\end{equation}
by (\ref{eq:rescaling}). In essence, the idea is to 
fix a function
$\f\in C^k_{c,r}(\R^n)$ such that $0\le\f\le1$,
 $\f\equiv 1$ on the unit ball, and to use
$E'(u_h)\cdot(\f^p u_h)=o(1)$
to reach a contradiction with (\ref{eq:contradiction}). However this can not be done if $p<k$,
 as $\f^p$ is not of class $C^k(\R^n)$. Thus we define
 $$
\Phi_\eps(x)=\left(\eps^2+\f(x)^2\right)^{p/2}-\eps^p~\!,
$$
where $\eps\in(0,1)$ is  fixed. Notice that $\Phi_\eps \in C^k_{c,r}(\R^n)$ and that $\Phi_\eps$ is a constant
in a neighborhood of $0$. Thanks to Lemma \ref{L:test1} we can compute
$E'(u_h)\cdot (\Phi_\eps u_h)$. Since 
the family 
$\Phi_\eps u_h$ is uniformly bounded in $\mathcal D^{k,\p}_{\rm r}(\R^n;|x|^\alpha dx)$
 as $h\to \infty$, then (\ref{eq:energy}) gives
\begin{equation}
\label{eq:testPhi}
e'_{k,p}(u)\cdot(\Phi_\eps u_h)=e'_{j,q}(u)\cdot(\Phi_\eps u_h)+o(1).
\end{equation}
We apply (\ref{eq:seconda}) in  Lemma \ref{L:test} to get
\begin{eqnarray*}
e'_{k,p}(u)\cdot(\Phi_\eps u_h)&=&
\irn|x|^\alpha|\nabla^k u_h|^{p-2}\nabla^k u_h\cdot\nabla^k(\Phi_\eps u_h)~\!dx\\
&=&
\irn|x|^\alpha|\nabla^k u_h|^{p}~\!\Phi_\eps~\! dx+o(1).
\end{eqnarray*}
From $\Phi_\eps\ge \f^p-\eps^p$, we infer
\begin{eqnarray*}
e'_{k,p}(u)\cdot(\Phi_\eps u_h)&\ge& \irn|x|^\alpha|\f\nabla^k u_h|^{p}~\! dx-c\eps^p+o(1)\\
&=&\irn|x|^\alpha|\nabla^k(\f u_h)|^{p}~\! dx-c\eps^p+o(1)
\end{eqnarray*}
by (\ref{eq:prima}), where $c=\sup_h\|u_h\|_{k,\alpha}$. Since $\f u_h\in\mathcal D^{k,\p}_{\rm r}(\R^n;|x|^\alpha dx)$
by Lemma \ref{L:test1}, from the definition of $\widetilde S$
we get
\begin{equation}
e'_{k,p}(u)\!\cdot\!(\Phi_\eps u_h)\ge 
\label{eq:left}
\tilde S\left(\irn|x|^{-\beta}|\nabla^j(\f u_h)|^{q}~\! dx\right)^{p/q}-c\eps^p+o(1).
\end{equation}
To estimate the right hand side of (\ref{eq:testPhi})
we  use (\ref{eq:seconda}) with $\alpha, k, p$ replaced by $\beta,j, q$,
respectively, to get
\begin{eqnarray*}
e'_{j,q}(u)\cdot(\Phi_\eps u_h)\!&=&\!
\irn|x|^{-\beta}|\nabla^j u_h|^{q-2}\nabla^j u_h\cdot\nabla^j(\Phi_\eps u_h)~\!dx\\
&=&\!
\irn|x|^{-\beta}|\nabla^j u_h|^{q-p}|\nabla^j u_h|^{p}~\!\Phi_\eps~\! dx+o(1).
\end{eqnarray*}
Thus
$$
e'_{j,q}(u)\!\cdot\!(\Phi_\eps u_h)\le 
\left(\int_{\{|x|<2\}}\!\!|x|^{-\beta}|\nabla^j \!u_h|^{q}~\!dx\right)^{\!\!\frac{q-p}{q}}\!\!\!
\left(\irn\!|x|^{-\beta}|\nabla^j \!u_h|^{q}(\Phi_\eps)^{\!\frac{q}{p}}~\!dx\right)^{\frac{p}{q}}\!\!\!+o(1)
$$
by H\"older inequality. From (\ref{eq:rescaling}), Lemma \ref{L:calculus} and (\ref{eq:prima}) we obtain
\begin{eqnarray}
\nonumber
e'_{j,q}(u)\cdot(\Phi_\eps u_h)&\le&
\frac{1}{2}~\!\tilde S~\!
\left(\irn\!|x|^{-\beta}|\nabla^j \!u_h|^{q}(\f^q+c\eps)~\!dx\right)^{\frac{p}{q}}\!\!\!+o(1)\\
\nonumber
&=& \frac{1}{2}~\!\tilde S~\!
\left(\irn\!|x|^{-\beta}|\nabla^j (\f u_h)|^{q}~\!dx+c\eps\right)^{p/q}\!\!+o(1)\\
&=& \frac{1}{2}~\!\tilde S~\!
\left(\irn\!|x|^{-\beta}|\nabla^j (\f u_h)|^{q}~\!dx\right)^{p/q}+c\eps^{p/q}+o(1).
\label{eq:right}
\end{eqnarray}
Comparing (\ref{eq:right}), (\ref{eq:testPhi}) and
(\ref{eq:left}) we conclude that
\begin{eqnarray*}
c\eps^{p/q}&+&\frac{1}{2}~\!\tilde S~\!\left(\irn\!|x|^{-\beta}|\nabla^j (\f u_h)|^{q}~\!dx\right)^{p/q}\\
&&\quad\quad\quad\ge e'_{j,q}(u)\cdot(\Phi_\eps u_h)+o(1)=
e'_{k,p}(u)\cdot(\Phi_\eps u_h)+o(1)\\
&&\quad\quad\quad\ge
\tilde S\left(\irn|x|^{-\beta}|\nabla^j(\f u_h)|^{q}~\! dx\right)^{p/q}-c\eps^p+o(1).
\end{eqnarray*}
Since $\eps^p=o(\eps^{p/q})$, we infer
$$
\irn|x|^{-\beta}|\nabla^j(\f u_h)|^{q}~\! dx\le 
c\eps+o(1),
$$
that together with (\ref{eq:contradiction}) implies
$$
\left(\frac{1}{2}~\!\widetilde S\right)^{\frac{q}{q-p}}\!\!+o(1)=
\int_{\{|x|<1\}}\!|x|^{-\beta}|\nabla^j u_h|^{q}~\!dx\le
\irn\!|x|^{-\beta}|\nabla^j (\f u_h)|^{q}~\!dx\le
c\eps+o(1),
$$
as $\f\equiv 1$ on the unit ball. The desired contradiction is achieved by choosing
$\eps>0$ small enough.
The proof is complete.
\QED

\begin{Remark}
\label{R:newSobolev}
Theorem \ref{T:extremalsk} includes $k$ Sobolev constants
without weights. Assume $n>kp$, take an integer $j\in\{0,,...,k-1\}$ and
$\alpha=0$. In addition, assume that $q$ equals the $(k-j)$th order
critical exponent
$$
p^{(k-j)*}:=\frac{np}{n-(k-j)\p}.
$$
Then $\beta_{k-j,q}=0$. Taking Remark \ref{R:densityk} into account, by
Theorem \ref{T:extremalsk} we have that
the radial Sobolev constant
$$
S_{k,p}^{{(k-j)*}}:=
\inf_{\scriptstyle u\in\mathcal D^{k,\p}_{\rm r}(\R^n)\atop\scriptstyle u\ne 0}
\frac{\displaystyle
\int_{\R^n}|\nabla^k u|^\p dx}
{\left(\displaystyle\int_{\R^n}|\nabla^ju|^{\frac{\p n}{n-(k-j)\p}}~\! dx\right)^{
\frac{n-(k-j)\p}{n}}}
$$
is positive and achieved in $\mathcal D^{k,\p}_{\rm r}(\R^n)$.
\end{Remark}

\appendix

\section{\!\!\!\!\!\!ppendix}
Here we prove some compactness and technical results that have been used in the
proof of Theorem \ref{T:extremalsk}. We always assume  that (\ref{eq:nondeg_gen}) is satisfied.

\begin{Lemma}
\label{L:test1}
Let $u\in \mathcal D^{k,\p}_{\rm r}(\R^n;|x|^\alpha dx)$ and $\Phi\in C^k_{c,r}(\R^n)$
such that $\Phi$ is  constant in a neighborhood of $0$. Then $\Phi u\in
\mathcal D^{k,\p}_{\rm r}(\R^n;|x|^\alpha dx)$ and 
$$
\irn|x|^\alpha|\nabla^k(\Phi u)|^p~\!dx\le c_\Phi \irn|x|^\alpha|\nabla^k u|^p~\!dx,
$$
where the constant $c_\Phi$ does not depend on $u$.
\end{Lemma}

\proof
Let $g=\mathcal T^{-1}_k u\in W^{k,p}(\R)$ and put $\tilde\Phi(s) :=\Phi(e^{-s})\in C^k(\R)$. 
Notice that
$\tilde\Phi(s)\equiv 0$ for
$s<<0$, $\tilde\Phi(s)$ is a constant for $s>>0$. Thus $\tilde\Phi g\in W^{k,p}(\R)$
and hence $\Phi u=\mathcal T_k(\tilde\Phi g)\in \mathcal D^{k,\p}_{\rm r}(\R^n;|x|^\alpha dx)$.
Finally, from Lemma \ref{L:RRSk} we infer that
$$
\irn|x|^\alpha|\nabla^k(\Phi u)|^p~\!dx\le c\|\tilde\Phi g\|_{W^{k,p}}^p
\le c_\Phi\|g\|_{W^{k,p}}^p\le c_\Phi \irn|x|^\alpha|\nabla^k u|^p~\!dx,
$$
and the Lemma is proved.
\QED

Next we need few results on weak convergence in 
$\mathcal D^{k,\p}_{\rm r}(\R^n;|x|^\alpha dx)$. Notice that 
$\mathcal D^{k,\p}_{\rm r}(\R^n;|x|^\alpha dx)$ is reflexive because 
it is topologically equivalent to $W^{k,p}(\R)$ (or because its
norm is uniformly convex). 

\begin{Lemma}
\label{L:Rellich}
Let $\Omega$ be a domain such that $\overline\Omega\subset\Rbuco$.
Then 
$\mathcal D^{k,\p}_{\rm r}(\R^n;|x|^\alpha dx)$ is compactly embedded into
 $W^{k-1,\tau}(\Omega)$ for any $\tau\ge p$.
\end{Lemma}

\proof
Use the Emden-Fowler transform and Rellich theorem for $W^{k,p}(I)$,
where $I$ is an appropriate bounded interval. 
\QED

\begin{Remark}
Lemma \ref{L:Rellich} is closely related to Theorem II.1 in \cite{PLL0}. Actually, the same argument
shows that 
$\mathcal D^{k,\p}_{\rm r}(\R^n;|x|^\alpha dx)$ is compactly embedded into
$C^{k-1}(\overline{\Omega})$. 
\end{Remark}

\begin{Lemma}
\label{L:test}
Let $\Phi\in C^k_{c,r}(\R^n)$ be a given function, 
such that $\Phi$ is  constant in a neighborhood of $0$.
If $u_h\weak 0$ in $\mathcal D^{k,\p}_{\rm r}(\R^n;|x|^\alpha dx)$
then
\begin{gather}
\label{eq:prima}
\irn|x|^\alpha|\nabla^k(\Phi u_h)-\Phi\nabla^k u_h|^p~dx=o(1)\\
\label{eq:seconda}
\irn|x|^\alpha|\nabla^k u_h|^{p-2}\nabla^k u_h\cdot\nabla^k(\Phi u_h)~\!dx=
\irn|x|^\alpha|\nabla^k u_h|^{p}~\!\Phi~\! dx+o(1)~\!.
\end{gather}
\end{Lemma}

\proof
Note that $\nabla^k(\Phi u_h)=\Phi\nabla^k u_h+\Psi D_{k-1}u_h$,
where $\Psi\in C^0_c(\Rbuco)$ and
$D_{k-1}$ is a differential operator of order $k-1$. Thus (\ref{eq:prima})
holds by Lemma \ref{L:Rellich}. To prove (\ref{eq:seconda}) use (\ref{eq:prima}) and
H\"older inequality.
\QED

We conclude this appendix with an elementary lemma.

\begin{Lemma}
\label{L:calculus}
Let $\eps,\f\in[0,1]$, $1<p<q$. Then there exists a constant $c\ge 0$ such that
$\left[\left(\eps^2+\f^2\right)^{p/2}-\eps^p\right]^{q/p}
\le\f^q+ c\eps$.
\end{Lemma}

\proof
Fix $\f$ and put $\Phi(\eps)=\left[\left(\eps^2+\f^2\right)^{p/2}-\eps^p\right]^{q/p}$. If $p\le2$ then
$\Phi$ is non increasing and hence the conclusion in the lemma holds with $c=0$. 
If $2<p<q$ it suffices to notice that $\Phi$ is differentiable at $\eps=0$.
\QED

\label{References}

\end{document}